\documentclass[leqno,twoside]{article}
\usepackage{amsmath,amssymb,amsthm}
\usepackage[margin=2.0cm]{geometry} % showframe,
\usepackage{titlesec,hyperref}
\usepackage{fancyhdr}

\usepackage{color}

\fancypagestyle{plain}
{
\fancyhf{}

}

\titleformat{\subsection}{\it}{\thesubsection.\enspace}{1pt}{}

\newtheorem{theo}{Theorem}[section]
\newtheorem{lemm}[theo]{Lemma}
\newtheorem{defi}[theo]{Definition}

\newtheorem{rema}[theo]{Remark}
\numberwithin{equation}{section}

\allowdisplaybreaks

\def\th2{\frac{\theta}{2}}

\begin{document}

\title{Global well-posedness for Euler-Nernst-Planck-Possion system \\in dimension two
\hspace{-4mm}
}

\author{Zeng Zhang$^1$
\quad Zhaoyang Yin$^2$ \\[10pt]
Department of Mathematics, Sun Yat-sen University,\\
510275, Guangzhou, P. R. China.\\[5pt]
}
\footnotetext[1]{Email: \it zhangzeng534534@163.com}
\footnotetext[2]{Email: \it mcsyzy@mail.sysu.com.cn}
\date{}
\maketitle
%\hrule
\begin{abstract}
In this paper, we study the Cauchy problem of the Euler-Nernst-Planck-Possion system. We obtain global well-posedness for the system in dimension $d=2$ for any initial data in $H^{s_1}(\mathbb{R}^2)\times H^{s_2}(\mathbb{R}^2)\times H^{s_2}(\mathbb{R}^2)$
under certain conditions of $s_1$ and $s_2$.

\vspace*{5pt}
%35K15Initial value problems for second-order parabolic equations
%35Q35 PDEs in connection with fuid mechanics

\vspace*{5pt}
\end{abstract}

\vspace*{10pt}
%\phantomsection
%\addcontentsline{toc}{section}{\contentsname}
%Ìí¼ÓĿ¼µ½ÊéÇ©
\tableofcontents

\section{Introduction}
In this paper, we study the Cauchy problem of the following nonlinear system:
\begin{align}\tag{01}\label{s1}
\left\{
\begin{array}{l}
u_t+u\cdot \nabla u-\nu \triangle u+\nabla P=\triangle\phi\nabla\phi, \quad t>0,\,x \in \mathbb{R}^d, \\[1ex]
\nabla\cdot u=0,\quad t>0,\,x \in \mathbb{R}^d,\\[1ex]
n_t+u\cdot \nabla n=\nabla\cdot(\nabla n-n\nabla\phi), \quad t>0,\,x \in \mathbb{R}^d,\\[1ex]
p_t+u\cdot \nabla p=\nabla\cdot(\nabla p+p\nabla\phi), \quad t>0,\,x \in \mathbb{R}^d,\\[1ex]
\triangle\phi=n-p,\quad t>0,\,x \in \mathbb{R}^d,\\[1ex]
(u,n,p)|_{t=0}=(u_0,n_0,p_0),\quad x \in \mathbb{R}^d.
\end{array}
\right.
\end{align}
Here $u(t,x)$ is a vector in $\mathbb{R}^d,$ $P(t,x),n(t,x),p(t,x)$ and $\phi(t,x)$ are scalars. The first two equations of the system (\ref{s1}) are the conservation equations of the incompressible flow. $u$ denotes the velocity filed, $P$ denotes the pressure, $\nu\geq 0$ denotes the fluid viscosity and $\phi$ denotes the electrostatic potential caused by the net charged particles. The third and the fourth equations of the system (\ref{s1}), which are the Nernst-Planck equations modified by the convective terms $u\cdot\triangle n$ and $u\cdot\nabla p,$ model the balance between diffusion and convective transport of charge densities by flow and electric fields. $n$ and $p$ are the densities of the negative and positive charged particles. They are coupled by the Poisson equation (the fifth equation).
The system (\ref{s1}) arises from electrohydrodynamics, which describing the dynamic coupling between incompressible flows and diffuse charge systems finds application in biology, chemistry and pharmacology. See \cite{Bazant,Joseph,Lin,Newman} for more details.

If the fluid viscosity $\nu>0,$ The above system (\ref{s1}) is the so called Navier-Stokes-Nernst-Planck-Possion ($NSNPP$) system, and it has been studied by several authors. Schmuck \cite{Schmuck} and Ryham \cite{Ryham} obtained the global existence of weak solutions in a bounded domain $\Omega$ in dimension $d\leq 3$ with Neumann and Dirichlet boundary conditions respectively. By using elaborate energy analysis, Li \cite{Li} studied the quasineutral limit in periodic domain. When $\Omega=\mathbb{R}^n,$ Joseph \cite{Joseph} established the existence of a unique smooth local solution for smooth initial dada by making using of Kato's semigroup ideas. The author also established the stability under the inviscid limit $\nu\rightarrow 0.$ Zhao et al. \cite{zhaoW,zhaoW3,zhaoG,zhaoW2} studied the local and global well-posedness in the critical Lebesgue spaces, modulation spaces, Triebel-Lizorkin spaces and Besov spaces by using the Banach fixed point theorem.

If,  on the other hand, $\nu=0,$ the above system (\ref{s1}) is the Euler-Nernst-Planck-Possion ($ENPP$) system. Recently, Zhang and Yin \cite{zy} proved the local well-posedness for the $ENPP$ system in Besov spaces in dimension $d\geq 2.$

The purpose of this paper is to get the global existence for the $ENPP$ system in dimension $d=2.$ Motivated by \cite{keben} for the study of the Euler system, we first introduce the following modified system
\begin{align}\tag{02}\label{s3}
\left\{
\begin{array}{l}
u_t+u\cdot \nabla u+\Pi(u,u)=\mathcal{P}\big((\nabla\cdot\xi)\xi\big),  \\[1ex]
n_t+\nabla\cdot (u n)-\triangle n=-\nabla\cdot(n\xi),\\[1ex]
p_t+\nabla\cdot (u p)-\triangle p=\nabla\cdot(p\xi), \\[1ex]
\xi=-\nabla(-\triangle)^{-1}(n-p),%,\\[1ex]
%(u,n,p,\xi)|_{t=0}=(u_0,n_0,p_0,\xi_0).
\end{array}
\right.
\end{align}
where $\mathcal{P}$ is the Leray projector defined as $\mathcal{P}=Id+\nabla (-\triangle)^{-1}\nabla\cdot,$ and $\Pi(\cdot,\cdot)$ is a bilinear operator defined by
$
  \Pi(u,v)=\sum_{j=1}^{5}\Pi_j(u,v),
$
with \begin{align*}&\Pi_1(u,v)=\nabla|D|^{-2}T_{\partial_iu^j} \partial_jv^i, &\Pi_2(u,v)&=\nabla|D|^{-2}T_{\partial_jv^i} \partial_iu^j,\\
&\Pi_3(u,v)=\nabla|D|^{-2}\partial_i\partial_j(I-\triangle_{-1})R(u^i,v^j),
&\Pi_4(u,v)&=\theta E_d\ast\nabla \partial_i\partial_j\triangle_{-1}R(u^i,v^j),\\
&\Pi_5(u,v)=\nabla \partial_i\partial_j\big((1-\theta ) E_d\big)\ast\triangle_{-1}R(u^i,v^j).\end{align*}
Here $\theta$ is a function of $\mathcal{D}(B(0,2))$ with value $1$ on $B(0,1),$ $E_d$ stands for the fundamental solution of $-\triangle,$ and $|D|^{-2}$ denotes the Fourier multiplier with symbol $|\xi|^{-2}.$ See Section 2 for the definitions of $T$ and $R.$ \\
We deduce form the second to the fourth equations of the system (\ref{s3}) that the dynamic equations of $(n+p, \xi)$ are
\begin{align*}
\left\{
\begin{array}{l}
(n+p)_t+\nabla\cdot \big(u (n+p)\big)-\triangle (n+p)=-\nabla\cdot\big((\nabla\cdot\xi)\xi\big), \\[1ex]
\xi_t-\triangle \xi+(-\nabla(-\triangle)^{-1}\nabla\cdot)\big(u(\nabla\cdot\xi)\big)=-(-\nabla(-\triangle)^{-1}\nabla\cdot)\big((n+p)\xi\big).
\end{array}
\right.
\end{align*}
Denote $\mathcal{L}=-\nabla (-\triangle)^{-1}\nabla\cdot=Id-\mathcal{P}.$ We then introduce the following system
\begin{align}\tag{03}\label{s2}
\left\{
\begin{array}{l}
u_t+u\cdot \nabla u+\Pi(u,u)=\mathcal{P}\big((\nabla\cdot\xi)\xi\big),  \\[1ex]
z_t+\nabla\cdot (u z)-\triangle z=-\nabla\cdot\big((\nabla\cdot\xi)\xi\big), \\[1ex]
\xi_t-\triangle\xi+\mathcal{L}\big(u(\nabla\cdot \xi))=-\mathcal{L}\big(z\xi).
\end{array}
\right.
\end{align}
Note that, by means of basic energy argument, the terms $\langle\mathcal{P}\big((\nabla\cdot\xi)\xi\big), u\rangle$ and $\langle\mathcal{L}\big(u(\nabla\cdot \xi)), \xi\rangle$ can be canceled out, which plays an important role in the proof of global existence.

We point out that in \cite{zy}, the term $\nabla\phi=\nabla(-\Delta)^{-1}(p-n)$ was controlled by $n-p$ through the Hardy-Littlewood-Sobolev inequality, i.e., $\|\nabla\phi\|_{L^q}\lesssim\|n-p\|_{L^p}$ with $1<p<d.$ Whereas, in this paper $n-p=\nabla\cdot\xi$ is controlled by $\xi$ through $\|n-p\|_{H^s}\lesssim\|\xi\|_{H^{s+1}}.$ Hence, these two papers solve the $ENPP$ system in different function spaces.

We can now state our main results:
\begin{theo}\label{a1}
Let $d\geq 2,~(s_1,s_2)\in\mathbb{R}^2,$ satisfying
\begin{align}\label{jibentiaojian}
  s_1>1+\frac{d}{2}, \textit{and}~  s_2+\frac{3}{2}> s_1\geq s_2+1.
\end{align}
There exists constants $c$ and $r$, depending only on $s_1,s_2$ and $d,$ such that for $(u_0,z_0,\xi_0)\in H^{s_1}(\mathbb{R}^d)\times H^{s_2}(\mathbb{R}^d)\times H^{s_2+1}(\mathbb{R}^d),$ with $\nabla\cdot u_0=0,$ $\xi_0=-\nabla(-\triangle)^{-1}a_0$ for some $a_0\in H^{s_2}(\mathbb{R}^d),$ and $z_0\pm\nabla \cdot \xi_0\geq 0$, there exists a time $$T\geq \frac{c}{1+(\|u_0\|_{ H^{s_1}(\mathbb{R}^d)}+\|z_0\|_{H^{s_2}(\mathbb{R}^d)}+\|\xi_0\|_{H^{s_2+1}(\mathbb{R}^d)})^r},$$ such that the system (\ref{s2}) has a unique solution $(u,z,\xi)$ on $[0,T]\times\mathbb{R}^d$ satisfying
$$(u,z,\xi)\in \widetilde{L}^\infty_{T}(H^{s_1}(\mathbb{R}^d))\times \big(\widetilde{L}^\infty_{T}(H^{s_2}(\mathbb{R}^d))
\cap\widetilde{L}^1_{T}(H^{s_2+2}(\mathbb{R}^d))\big)\times \big(\widetilde{L}^\infty_{T}(H^{s_2+1}(\mathbb{R}^d))
\cap\widetilde{L}^1_{T}(H^{s_2+3}(\mathbb{R}^d))\big),$$
 and $(u,z,\xi)$ is continuous in time with values in $H^{s_1}\times H^{s_2}\times H^{s_2+1}.$\\
 Moreover, $\nabla\cdot u=0,$ $z\pm\nabla\cdot \xi\geq 0,~a.e.~on ~[0,T]\times \mathbb{R}^d,$ and $\mathcal{L}\xi=\xi.$\\
Finally, if $d=2$ and $s_2>1,$ then the solution $(u,z,\xi)$ is global.\end{theo}
\begin{rema}
We mention that the restriction $s_1>1+\frac{d}{2}$ is due to the same reasons as illustrated for the Euler equation in \cite{keben}. $s_2+\frac{3}{2}> s_1\geq s_2+1$ is caused by the coupling between $u$ and $\xi.$ In fact, owing to the properties of the transport flow, $(u,\xi)$ is expected to be in $ \widetilde{L}^\infty_{T}(H^{s_1}(\mathbb{R}^d))\times \widetilde{L}^r_{T}(H^{s_2+1+\frac{2}{r}}(\mathbb{R}^d))$ with $r\in [1,\infty].$ Due to the product laws in Besov spaces, $ab$ is less regular than $a$ or $b$. Thus in order to control the term $(\nabla\cdot\xi)\xi$ in the first equation of the system, we have to assume
$s_2+1+\frac{2}{r_1}\geq (>)s_1,~s_2+\frac{2}{r_1'}> (\geq) s_1,$
for some $r_1\in [1,\infty]$ and $\frac{1}{r_1}+\frac{1}{r_1'}=1,$ which implies that $s_2+\frac{3}{2}>s_1.$ Similar reason for the term $u(\nabla\cdot\xi)$ requires $s_1\geq s_2+1.$
\end{rema}
\begin{theo}\label{a3}
Let $d= 2,~(s_1,s_2)\in\mathbb{R}^2,$
\begin{align}\label{xinde}
  s_1>2, ~s_2>1,\textit{and}~  s_2+\frac{3}{2}> s_1\geq s_2+1.
\end{align}
Then for any $(u_0,n_0,p_0)\in H^{s_1}(\mathbb{R}^2)\times H^{s_2}(\mathbb{R}^2)\times H^{s_2}(\mathbb{R}^2),$ with $\nabla\cdot u_0=0,$ $\nabla(-\triangle)^{-1}(n_0-p_0)\in H^{s_2+1}(\mathbb{R}^2),$ and $n_0,p_0\geq 0$, the $ENPP$ system has a solution $(u,n,p,P,\phi)$ on $\mathbb{R}^+\times\mathbb{R}^2$ satisfying
\begin{align*}&(u,n,p)\in \widetilde{L}^\infty(\mathbb{R}^+;H^{s_1}(\mathbb{R}^2))\times \Big(\widetilde{L}^\infty(\mathbb{R}^+;H^{s_2}(\mathbb{R}^2))
\cap\widetilde{L}^1(\mathbb{R}^+;H^{s_2+2}(\mathbb{R}^2))\Big)^2,\\
&-\nabla(-\triangle)^{-1}(n-p)\in \widetilde{L}^\infty(\mathbb{R}^+;H^{s_2+1}(\mathbb{R}^2))
\cap\widetilde{L}^1(\mathbb{R}^+;H^{s_2+3}(\mathbb{R}^2)),\\
&P,\Phi\in L^\infty(\mathbb{R}^+;BMO(\mathbb{R}^2)).\end{align*}
Morever, if $(\widetilde{u},\widetilde{n},\widetilde{p},\widetilde{P},\widetilde{\phi})$ also satisfies the $ENPP$ system with the same initial data and belongs to the above class, then $(u,n,p)=(\widetilde{u},\widetilde{n},\widetilde{p}),$ and $(\nabla P,\nabla \phi)=(\nabla\widetilde{P},\nabla\widetilde{\phi}).$\\
Finally,
 $(u,n,p)$ is continuous in time with values in $H^{s_1}\times H^{s_2}\times H^{s_2},$
and $n.p\geq 0,~a.e.~on ~\mathbb{R}^+\times \mathbb{R}^2$.\end{theo}
\begin{rema}
  We mention that under an improved condition \ref{xinde}, Theorem \ref{a3} may hold true for the $NSNPP$ system. We will present this result in another paper.
\end{rema}
Throughout the paper, $C>0$ stands for a generic constant and $c>0$ a small constant. We shall sometimes use the notation $A\lesssim B$ to denote the relation $A \leq CB.$ For simplicity, we write $L^p$, $H^s$ and $B^s_{p,r}$ for the spaces $L^p(\mathbb{R}^d),~H^s(\mathbb{R}^d),$ and $B^s_{p,r}(\mathbb{R}^d)$, respectively.

The remain part of this paper is organized as follows.
In Section 2, we recall some basic facts about Littlewod-Paley theory and Besov spaces. Section 3 is devoted to the proof of Theorem \ref{a1}. Finally, we give a proof of Theorem \ref{a3} by using Theorem \ref{a1}.
\section{Preliminaries}
  \subsection{The nonhomogeneous Besov spaces }
We first define the Littlewood-Paley decomposition.
 \begin{lemm}\cite{keben}
Let $\mathcal{C}=\{\xi\in{\mathbb{R}^2}, ~\frac{3}{4}\leq|\xi|\leq\frac{8}{3} \}$ be an annulus. There exist radial functions $\chi$ and $\varphi$ valued in the interval $[0,1]$, belonging respectively to $\mathcal{D}(B(0,\frac{4}{3}))$ and $\mathcal{D}(\mathcal{C})$,  such that
\begin{align*}
\forall \xi \in{\mathbb{R}}^d,~~\chi(\xi)+\sum_{j\geq 0}\varphi(2^{-j}\xi)=1.\end{align*}\end{lemm}
The nonhomogeneous dyadic blocks $\triangle_j$ and the nonhomogeneous low-frequency cut-off operator $S_j$ are then defined as follows:
\begin{align*}&\triangle_ju=0~~\textit{if}~j\leq-2,~~~~~~~~~~~~~\triangle_{-1}u=\chi(D)u,\\
&\triangle_ju=\varphi(2^{-j}D)u~~\textit{if}~j\geq0,~~~~S_ju=\sum_{j'\leq j-1}\triangle_{j'}u,~~\textit{for}~j\in \mathbb{Z}.\end{align*}

We may now introduce the nonhomogeneous Besov spaces.
\begin{defi}\label{dingyi}
Let $s\in \mathbb{R}$ and $(p,r)\in[1,\infty]^2$. The nonhomogeneous Besov space $B^s_{p,r}$
consists of all tempered distributions $u$ such that
$$\|u\|_{B^s_{p,r}}\overset{def}{=}\Big\|(2^{js}\|\triangle_ju\|_{L^p})_{j\in \mathbb{Z}}\Big\|_{l^r(\mathbb{Z})}
<\infty.$$
\end{defi}
The Sobolev space can be defined as follows:
 \begin{defi}
 For $s\in \mathbb{R}$,
 \begin{align*}
H^s =\{u \in \mathcal{S}'; \|u\|_{H^s}=\Big(\sum_{j=-1}^{\infty}2^{2js}\|\triangle_ju\|_{L^2}^2\Big)^{\frac{1}{2}}<\infty\}.
\end{align*}
\end{defi}
\begin{rema}
  For any $s\in \mathbb{R}$, the Besov space $B^{s}_{2,2}$
coincides with the Sobolev space $H^s.$
\end{rema}
\begin{lemm}\label{Fadou}
The set $B^s_{p,r}$ is a Banach space, and satisfies the Fatou property, namely, if $(u_n)_{n\in N}$ is a
bounded sequence of $B^s_{p,r}$, then an element $u$ of $B^s_{p,r}$ and a subsequence $u_{\psi(n)}$ exist such that
$$\underset{n\rightarrow\infty}{\lim}~u_{\psi(n)}=u~~in~~\mathcal{S}'~~~and ~~~\|u\|_{B^s_{p,r}}\leq C \underset{n\rightarrow\infty}{\liminf} \|u_{\psi(n)}\|_{B^s_{p,r}}.$$\end{lemm}
In addition to the general time-space $L^{\rho}_T(B^s_{p,r})$, we introduce the following mixed time-space $\widetilde{L}^{\rho}_T(B^s_{p,r}).$
\begin{defi}
For all $T>0,~s\in\mathbb{R},$ and $1\leq r,\rho\leq\infty$, we define the space $\widetilde{L}^{\rho}_T(B^s_{p,r})$ the set of tempered distributions $u$ over $(0,T)\times \mathbb{R}^d,$ such that
$$\|u\|_{\widetilde{L}^{\rho}_T(B^s_{p,r})}\overset{def}{=}\|2^{js}\|\triangle_ju\|_{L^{\rho}_T(L^p)}\|
_{l^r(\mathbb{Z})}<\infty.$$
\end{defi}
\noindent It follows from the Minkowski inequality that
\begin{align*}
  \|u\|_{L^\rho_T({B}_{p,r}^{s})}\leq \|u\|_{\widetilde{L}^\rho_T({B}_{p,r}^{s})}~if~r\leq \rho,~~\|u\|_{\widetilde{L}_T^\rho({B}_{p,r}^{s})}\leq \|u\|_{L^\rho_T({B}_{p,r}^{s})}~if~r\geq \rho.
\end{align*}

  Let's then recall Bernstein-Type lemmas.
\begin{lemm}\label{Bi}\cite{keben} (Bernstein inequalities)
Let $\mathcal{C}$ be an annulus and $\mathcal{B}$ a ball. A constant $C$ exists such that for any nonnegative integer $k$, any couple $(p,q)$ in $[1,\infty]^2$ with $q\geq p\geq1$, and any function u of $L^p$, we have
\begin{align*}
&Supp\, \widehat{u}\subset \lambda\mathcal{B}\Rightarrow\,\underset{|\alpha|=k}{\sup}\,\|\partial^{\alpha}u\|_{L^q}
\leq C^{k+1}\lambda^{k+d(\frac{1}{p}-\frac{1}{q})}\|u\|_{L^p},\\
&Supp\,\widehat{u}\subset \lambda\mathcal{C}\Rightarrow\,C^{-k-1}\lambda^k\|u\|_{L^p}\leq\underset{|\alpha|=k}{\sup}\,\|\partial^{\alpha}u\|_{L^q}
\leq C^{k+1}\lambda^k\|u\|_{L^p}.\end{align*}
\end{lemm}
%\begin{lemm}\label{DL}\cite{DanchinL} (A Bernstein-like inequality)
%Let $1< p <\infty$ and $u\in L^p(\mathbb{R}^d),$ such that $Supp
%\,
%\widehat{u}\in C(0,R_1,R_2)~(with~ 0 < R_1 < R_2).$ There exists a constant c depending
%only on $d$ and $R_2/R_1$, such that
%\begin{align*}
%c\frac{R_1^2}{p^2}\int_{\mathbb{R}^d}|u|^pdx\leq\int_{\mathbb{R}^d} |\nabla u|^2|u|^{p-2}dx=-\frac{1}{p-1}\int_{\mathbb{R}^d}
%\triangle u|u|^{p-2}udx.
%\end{align*}
%\end{lemm}

 We state the following embedding and interpolation inequalities.
\begin{lemm}\cite{keben}
  Let $1\leq p_1\leq p_2\leq\infty$ and $\leq r_1\leq r_2\leq\infty.$ Then for any real number $s,$ we have
  $${B}^s_{p_1,r_1}\hookrightarrow{B}^{s-d(\frac{1}{p_1}-\frac{1}{p_2})}_{p_2,r_2}.$$
  \end{lemm}
  \begin{lemm}\cite{keben}
  If $s_1$ and $s_2$ are real numbers such that $s_1<s_2,$  $\theta\in(0,1)$ and $~1\leq p,r\leq\infty,$ then we have
 \begin{align*}
 B^{s_2}_{p,\infty}\hookrightarrow B^{s_1}_{p,1},~~~\textit{and}~~~
 \|u\|_{{B}^{\theta s_1+(1-\theta)s_2}_{p,r}}\leq \|u\|_{B^{s_1}_{p,r}}^\theta\|u\|_{B^{s_2}_{p,r}}^{1-\theta}.
  \end{align*}
  \end{lemm}

In the sequel, we will frequently use the Bony decomposition:
$$uv=T_vu+T_uv+R(u,v),$$
with
\begin{align*}&R(u,v)=\underset{|k-j|\leq1}{\sum}\triangle_ku\triangle_jv,\\
&T_uv=\underset{j\in \mathbb{Z}}{\sum}S_{j-1}u\triangle_jv
=\underset{j\geq1}{\sum}S_{j-1}u\triangle_j\big((Id-\triangle_{-1})v\big)~~
,\end{align*}
where operator $T$ is called ``paraproduct", whereas $R$ is called ``remainder".
\begin{lemm}\label{T}
A constant $C$ exists which satisfies the following inequalities for any couple of real numbers $(s,t)$ with t negative and any $(p,p_1,p_2,r,r_1,r_2)$ in $[1,\infty]^6$:
\begin{align*}&\|T\|_{\mathcal{L}(L^{p_1}\times {B}^s_{p_2,r};{B}^s_{p,r})}\leq C^{|s|+1},\\
&\|T\|_{\mathcal{L}({B}^t_{p_1,r_1}\times {B}^s_{p_2,r_2};{B}^{s+t}_{p,r})}\leq \frac{C^{|s+t|+1}}{-t},
\end{align*}
with $\frac{1}{p}\overset{def}{=}\frac{1}{p_1}+\frac{1}{p_2}\leq1,~ \frac{1}{r}\overset{def}{=}min\{1,\frac{1}{r_1}+\frac{1}{r_2}\}.$
\end{lemm}
\noindent{Proof.} The proof of this lemma can be easily deduced from substituting the estimate
$$\|S_{j-1}u\triangle_j v\|_{L^p}\leq \|S_{j-1}u\|_{L^{p_1}}\|\triangle_j v\|_{L^{p_2}},$$
for the estimate
$$\|S_{j-1}u\triangle_j v\|_{L^p}\leq \|S_{j-1}u\|_{L^{\infty}}\|\triangle_j v\|_{L^{p}}$$
in the proof of Theorem 2.82 in \cite{keben}. It is thus omitted.\qed
\begin{lemm}\label{R}\cite{keben}
A constant $C$ exists which satisfies the following inequalities. Let $(s_1,s_2)$ be in $\mathbb{R}^2$ and
$(p_1,p_2,r_1,r_2)$ be in $[1,\infty]^4$. Assume that
$$\frac{1}{p}\overset{def}{=}\frac{1}{p_1}+\frac{1}{p_2}\leq1~~and~~
\frac{1}{r}\overset{def}{=}\frac{1}{r_1}+\frac{1}{r_2}\leq1.$$
If $s_1+s_2>0$, then we have, for any $(u,v)$ in ${B}^{s_1}_{p_1,r_1}\times {B}^{s_2}_{p_2,r_2}$,
$$\|R(u,v)\|_{{B}^{s_1+s_2}_{p,r}}\leq\frac{C^{|s_1+s_2|+1}}{s_1+s_2}\|u\|_{{B}^{s_1}_{p_1,r_1}}\|v\|_{{B}^{s_2}_{p_2,r_2}}.$$
If $r=1$ and $s_1+s_2=0$, then we have, for any $(u,v)$ in ${B}^{s_1}_{p_1,r_1}\times {B}^{s_2}_{p_2,r_2}$,
$$\|R(u,v)\|_{{B}^0_{p,\infty}}\leq C\|u\|_{{B}^{s_1}_{p_1,r_1}}\|v\|_{{B}^{s_2}_{p_2,r_2}}.$$
\end{lemm}
\begin{lemm}\label{s_2+1}
 Let $s+\frac{1}{2}>\frac{d}{2}.$ A constant $C$ exists such that
 \begin{align*}
   &\|uv\|_{H^s}\lesssim\|u\|_{H^{s+\frac{1}{2}}}\|v\|_{H^s},\\
   &\|uv\|_{H^{s+1}}\lesssim\|u\|_{H^{s+1}}\|v\|_{H^{s+1}}.
 \end{align*}
\end{lemm}
{\noindent Proof.} By using Bony's decomposition combined with Lemmas \ref{T}-\ref{R}, we have
\begin{align*}
  \|uv\|_{H^{s}}\lesssim&\|T_uv\|_{H^{s}}+\|R(u,v)\|_{H^{s}}+\|T_vu\|_{H^{s}}\\
  \lesssim&\|u\|_{L^\infty}\|v\|_{H^{s}}+\|u\|_{B^0_{\infty,\infty}}\|v\|_{H^{s}}
  +\|v\|_{B^{-\frac{1}{2}}_{\infty,\infty}}\|u\|_{H^{s+\frac{1}{2}}}\\\lesssim&\|u\|_{L^\infty}\|v\|_{H^{s}}
  +\|v\|_{H^{\frac{d}{2}-\frac{1}{2}}}\|u\|_{H^{s+\frac{1}{2}}}\\\lesssim&
  \|u\|_{H^{s+\frac{1}{2}}}\|v\|_{H^{s}},
\end{align*}
where we have used $s>\frac{d}{2}-\frac{1}{2}>0,$ and $H^{s+\frac{1}{2}}\hookrightarrow L^\infty.$ Similarly,
\begin{align*}
  \|uv\|_{H^{s+1}}
  \lesssim\|u\|_{L^\infty}\|v\|_{H^{s+1}}
  +\|v\|_{L^\infty}\|u\|_{H^{s+1}}\lesssim
  \|u\|_{H^{s+1}}\|v\|_{H^{s+1}},
\end{align*}
where we have used $H^{s+1}\hookrightarrow L^\infty.$ We thus obtain the desired inequalities. \qed\\
We mention that all the properties of continuity for the paraproduct and remainder remain true in the mixed time-space $\widetilde{L}^{\rho}_T(B^s_{p,r}).$

Finally, we state the following commutator estimates.
\begin{lemm}\label{jiaohuan}\cite{keben}
 Let $v$ be a vector filed over $\mathbb{R}^d,$ define $R_j=[v\cdot\nabla,\triangle_j]f.$ Let $\sigma>0~(\textit{or}~\sigma>-1,$ if $ \nabla\cdot v=0),$ $1\leq r\leq\infty,$ $1\leq p\leq p_1\leq\infty,$ and $\frac{1}{p_2}=\frac{1}{p}-\frac{1}{p_1}.$ Then
 \begin{align*}
   \Big\|2^{j\sigma}\|R_j\|_{L^P}\Big\|_{l^r}\leq C\Big(\|\nabla v\|_{L^\infty}\|f\|_{B^{\sigma}_{p,r}}+\|\nabla f\|_{L^{p_2}}\|\nabla v\|_{B^{\sigma-1}_{p_1,r}}\Big).
 \end{align*}
\end{lemm}
\subsection{A priori estimates for transport and transport-diffusion equations}
Let us state some classical a priori estimates for transport equations and transport-diffusion equations.
\begin{lemm}\label{ts}\cite{keben}
Let $1\leq p\leq p_1\leq\infty,~1\leq r\leq\infty$. Assume that
\begin{align}\label{tiaojian}
s\geq-d\,min\left(\frac{1}{p_1},\frac{1}{p'}\right) \quad \textit{or} \quad s\geq-1-d\,min\left(\frac{1}{p_1},\frac{1}{p'}\right)~\textit{if}~\nabla\cdot v=0
\end{align}
with strict inequality if $r<\infty$.

There exists a constant $C$, depending only on $d, p, p_1, r$ and $s$, such that for all solutions
$f\in L^{\infty}([0,T];B^s_{p,r})$ of the transport equation
\begin{align}
\left\{
\begin{array}{l}
\partial_tf+v\cdot\nabla f=g\\
f_{|t=0}=f_0,
\end{array}
\right.
\end{align}
with initial data $f_0$ in $B^s_{p,r}$, and $g$ in $L^1([0,T];B^s_{p,r})$, we have, for $a.e.\,t\in[0,T]$,
\begin{align}\label{,}\|f\|_{\widetilde{L}_t^{\infty}(B^s_{p,r})}\leq\left(\|f_0\|_{B^s_{p,r}}+
\int_0^t exp(-CV_{p_1}(t))\|g(t')\|_{B^s_{p,r}}dt'\right)exp(CV_{p_1}(t)),\end{align}
%or more accurately,
%\begin{align}\label{,,}\|f\|_{\widetilde{L}_t^{\infty}(B^s_{p,r})}\leq\|f_0\|_{B^s_{p,r}}+
%C\int_0^t V_{p_1}(t')\|f(t')\|_{\widetilde{L}_{t'}^{\infty}(B^s_{p,r})}dt'+
%\|g\|_{\widetilde{L}_t^{1}(B^s_{p,r})},\end{align}
with, if the inequality is strict in (\ref{tiaojian}),
\begin{align}
V'_{p_1}(t)=\left\{\begin{array}{l}\|\nabla v(t)\|_{B^{s-1}_{p_1,r}},~if~s>1+\frac{d}{p_1}~or~s=1+\frac{d}{p_1},~r=1,\\
\|\nabla v(t)\|_{B^{\frac{d}{p_1}}_{p_1,\infty}\cap L^{\infty}},~if~s<1+\frac{d}{p_1}
\end{array}\right.
\end{align}
and, if equality holds in (\ref{tiaojian}) and $r=\infty$,
$$V'_{p_1}=\|\nabla v(t)\|_{B^{\frac{d}{p_1}}_{p_1,1}}.$$
If $f=v$, then for all $s>0$ $(s>-1,$ if $\nabla\cdot u=0)$, the estimate (\ref{,}) %and (\ref{,,})
holds with
$$V'_{p_1}(t)=\|\nabla u\|_{L^{\infty}}.$$
%Moreover, $f\in C([0,T];B^s_{p,r}) ~if ~r<\infty;$\\
%$~~~~~~~~~~~~~f\in C([0,T];B^{s'_1}_{p,r}) ~if ~r=\infty,~for ~any ~s'_1<s.$
\end{lemm}
\begin{lemm}\label{tds}\cite{keben}
Let $1\leq p_1\leq p\leq\infty,~1\leq r\leq\infty,~s\in\mathbb{R}$ satisfy (2.10), and let $V_{p_1}$ be defined as in Lemma \ref{ts}.

There exists a constant $C$ which depends only on $d, r, s$ and $s-1-\frac{d}{p_1}$ and is such that for any smooth solution $f$ of the transport diffusion equation
\begin{align}
\left\{
\begin{array}{l}
\partial_tf+v\cdot\nabla f-\nu\triangle f=g\\
f_{|t=0}=f_0,
\end{array}
\right.
\end{align} we have
\begin{align*}\nu^{\frac{1}{\rho}}\|f\|_{\widetilde{L}^{\rho}_T(B^{s+\frac{2}{\rho}}_{p,r})}\leq Ce^{C(1+\nu T)^{\frac{1}{\rho}}V_{p_1}(T)}\Big(&(1+\nu T)^{\frac{1}{\rho}}\|f_0\|_{B^s_{p,r}}\\
+&(1+\nu T)^{1+\frac{1}{\rho}-\frac{1}{\rho_1}}\nu^{\frac{1}{\rho_1}-1}\|g\|_
{\widetilde{L}^{\rho_1}_T(B^{s-2+\frac{2}{\rho_1}}_{p,r})}\Big),
\end{align*} where $1\leq\rho_1\leq\rho\leq\infty.$
\end{lemm}

\subsection{The operator $\Pi(\cdot,\cdot)$}
We recall some basic results for $\Pi(\cdot,\cdot).$ See \cite{keben} (Pages 296-300) for further details.
\begin{lemm}\label{pi}\cite{keben}
  For all $s>-1,$ and $1\leq p,r\leq\infty,$ there exists a constant C such that
  \begin{align*}
    \|\Pi(v,w)\|_{B^s_{p,r}}\leq C(\|v\|_{C^{0,1}}\|w\|_{B^s_{p,r}}+\|w\|_{C^{0,1}}\|v\|_{B^s_{p,r}}).
  \end{align*}
  Moveover, there exists a bilinear operator $P_\Pi$ such that
   $\Pi(v,w)=\nabla P_\Pi(v,w),$ and
   \begin{align*}
    \|P_\Pi(v,w)\|_{B^{s+1}_{p,r}}\leq C\Big(\|v\|_{C^{0,1}}\|w\|_{B^s_{p,r}}+\|w\|_{C^{0,1}}\|v\|_{B^s_{p,r}}\Big),~if ~ 1<p<\infty
%    \|P_\Pi(v,w)\|_{B^{s+1}_{p,r}+L^\infty_L}\leq C\Big(\|v\|_{C^{0,1}}\|w\|_{B^s_{p,r}}+\|w\|_{C^{0,1}}\|v\|_{B^s_{p,r}}\Big),~if ~ p=\infty.
  .\end{align*}
\end{lemm}
\begin{lemm}\label{deng}\cite{keben}
  For all $-1<s<\frac{d}{p}+1,$ and $1\leq p,r\leq\infty,$ we have
  \begin{align*}
    \|\Pi(v,w)\|_{B^s_{p,r}}\leq C\Big(\|v\|_{C^{0,1}}\|w\|_{B^s_{p,r}}+\|w\|_{B^{s-\frac{d}{p}}_{\infty,\infty}}\|\nabla v\|_
    {B^{\frac{d}{p}}_{p,r}}\Big).
  \end{align*}
  \end{lemm}
  \begin{lemm}\label{yyy}\cite{keben}
  For all $s>1,$ and $1\leq p,r\leq\infty,$ there exists a constant C such that
  \begin{align*}
    \|\nabla\cdot\Pi(v,w)+tr(Dv,Dw)\|_{B^{s-1}_{p,r}}\leq C\Big(\|\nabla\cdot v\|_{B^{0}_{\infty,\infty}}\|w\|_{B^s_{p,r}}+
    \|\nabla\cdot w\|_{B^{0}_{\infty,\infty}}\|v\|_{B^s_{p,r}}\Big).
   % In the limit case $s=1$ we have
  \end{align*}
  \end{lemm}
 \subsection{The space $B^1_{\infty,\infty}$ }
  The space $B^1_{\infty,\infty}$ plays an important role in dealing with the global existence. In this section, we introduce an interpolation inequality involving $B^1_{\infty,\infty}$.
 % \begin{defi}
%  The space $LL$ consists of those bounded functions $f$ such that
%  \begin{align*}
%    \|f\|_{LL}\overset{def}{=}\underset{0<|x-x'|\leq 1}{\sup}\frac{|f(x)-f(x')|}{|x-x'|(1-log|x-x'|)}<\infty.
%  \end{align*}
%\end{defi}
\begin{defi}
Let $\alpha$ be in $(0,1].$ A modulus of continuity is any nondecreasing
nonzero continuous function $\mu:[0,\alpha]\rightarrow \mathbb{R}^+$ such that $\mu(0) = 0.$ The
modulus of continuity $\mu$ is admissible if, in addition, the function $\Gamma$ defined
for $y \geq \frac{1}{\alpha}$ by
$$\Gamma(y) \overset{def}{=} y\mu(\frac{1}{y})$$
is nondecreasing and satisfies, for some constant $C$ and all $x \geq \frac{1}{\alpha},$
$$\int_x^\infty\frac{1}{y^2}\Gamma(y)dy\leq C\frac{\Gamma(x)}{x}. $$
\end{defi}
\begin{defi}Let $\mu$ be a modulus of continuity and $(X, d)$ a metric space.
We denote by $C_{\mu}(X)$ the set of bounded, continuous, real-valued functions $u$
over $X$ such that
\begin{align*}
    \|u\|_{C_{\mu}}\overset{def}{=}\|u\|_{L^\infty}+\underset{0<d(x,y)\leq \alpha}{\sup}\frac{|u(x)-u(x')|}{\mu(d(x,y))}<\infty.
  \end{align*}
  \end{defi}
\begin{rema}\cite{keben}
%1. $B^1_{\infty,\infty}$ is a subspace of the space $LL$ of log-Lipschitz functions. More precisely, $B^1_{\infty,\infty}\hookrightarrow LL.$ \\
%2. If $\mu(r) = r(1-log r)$, then $C_{\mu}(X)$ is the space
%$LL(X)$ of log-Lipschitz functions on $X$.
Let $\alpha=1.$ The function $\mu(r) = r(1-log r)$ is an admissible modulus of continuity, and the space $C_{\mu}$ contains $B^1_{\infty,\infty},$ more precisely, $B^1_{\infty,\infty}\hookrightarrow C_{\mu}.$
\end{rema}
\begin{lemm}\label{LL}\cite{keben}
  Let $\mu$ be an admissible modulus of continuity. There exists
a constant C such that for any $\varepsilon\in(0,1],$ $u$ in $C^{1,\varepsilon},$ and positive $\Lambda$, we
have
\begin{align*}
  \|\nabla u\|_{L^\infty}\leq C\
  \left(\frac{\|u\|_{C_{\mu}}+\Lambda}{\varepsilon}+\|u\|_{C_\mu}
  \Gamma\Big(\big(\frac{\|\nabla u\|_{C^{0,\varepsilon}}}{\|u\|_{C_{\mu}}+\Lambda}\big)^{\frac{1}{\varepsilon}}\Big)\right)
\end{align*}
whenever $\|u\|_{C_{\mu}}+\Lambda\leq (\frac{\alpha}{2})^{\varepsilon}\|\nabla u\|_{C^{0,\varepsilon}}.$
\end{lemm}
\section{Proof of Theorem \ref{a1}}
\hspace{0.5cm}To begin, we denote
$
\varepsilon=s_2+\frac{3}{2}-s_1,
~\textit{and}~ \varepsilon_0=\min(\frac{1}{2},\varepsilon).$
We mention that the condition $(\ref{jibentiaojian})$ implies that
\begin{align}\label{jingchangyong}
s_2+\frac{1}{2}>s_1-1>\frac{d}{2},
\end{align}
which will be frequently used.
\subsection{Existence for the system (\ref{s2})}
  \subsubsection{First step: Construction of approximate solutions and uniform bounds}
  In order to define a sequence $(u^m,z^m,\xi^m)|_{m\in\mathbb{N}}$ of global approximate solutions to the system (\ref{s2}), we use an iterative scheme.
  First we set $u^0=u_0,~z^0=e^{t\triangle}z_0,~\xi^0=e^{t\triangle}\xi_0.$ Thanks to Lemma \ref{tds}, it is easy to see that
    $$(u^0,z^0,\xi^0)\in \widetilde{L}^\infty_{loc}(\mathbb{R}^+;H^{s_1})\times \Big(\widetilde{L}^\infty_{loc}(\mathbb{R}^+;H^{s_2})\cap\widetilde{L}^1_{loc}(\mathbb{R}^+;H^{s_2+1})\Big)^2,$$
    and
   \begin{align*}
   &\|u^0\|_{\widetilde{L}^\infty_t(H^{s_1})}
   +\|z^0\|_{\widetilde{L}^\infty_t(H^{s_2})\cap\widetilde{L}^1_t(H^{s_2+2})}
   +\|\xi^0\|_{\widetilde{L}^\infty_t(H^{s_2+1})\cap\widetilde{L}^1_t(H^{s_2+3})}\\\leq &C(1+t)(\|u_0\|_{H^{s_1}}+\|z_0\|_{H^{s_2}}+\|\xi_0\|_{H^{s_2+1}}).
   \end{align*}
   Then, assuming that $$(u^m,z^m,\xi^m)\in \widetilde{L}^\infty_{loc}(\mathbb{R}^+;H^{s_1})\times \big(\widetilde{L}^\infty_{loc}(\mathbb{R}^+;H^{s_2})\cap\widetilde{L}^1_{loc}(\mathbb{R}^+;H^{s_2+2})\big)\times \big(\widetilde{L}^\infty_{loc}(\mathbb{R}^+;H^{s_2+1})\cap\widetilde{L}^1_{loc}(\mathbb{R}^+;H^{s_2+3})\big),$$
   we solve the following linear system:
  \begin{align}
\left\{
\begin{array}{l}
u^{m+1}_t+u^{m}\cdot \nabla u^{m+1}=-\Pi(u^{m},u^{m})-\mathcal{P}\big((\nabla\cdot\xi^{m})\xi^{m}\big),  \\[1ex]
z^{m+1}_t-\triangle z^{m+1}=-\nabla\cdot(u^{m} z^{m})-\nabla\cdot\big((\nabla\cdot\xi^{m})\xi^{m}\big), , \\[1ex]
\xi^{m+1}_t-\triangle \xi^{m+1}=-\mathcal{L}\big(u^{m}(\nabla\cdot\xi^{m})\big)-\mathcal{L}\big(z^{m}\xi^{m}\big), \\[1ex]
(u^{m+1},z^{m+1},\xi^{m+1})|_{t=0}=(u_0,z_0,\xi_0).
\end{array}
\right.
\end{align}
Lemma \ref{ts} ensures that
\begin{align}\label{u'}
  \|u^{m+1}\|_{\widetilde{L}^\infty_t(H^{s_1})}\lesssim & exp(C\int_0^t\|u^{m}\|_{H^{s_1}}dt')\Big(
  \|u_0\|_{H^{s_1}}\\\nonumber&+\|\Pi(u^{m},u^{m})\|_{\widetilde{L}^1_t(H^{s_1})}+
  \|\mathcal{P}\big((\nabla\cdot\xi^{m})\xi^{m}\big)\|_{\widetilde{L}^1_t(H^{s_1})}\Big).
\end{align}
 Using Lemma \ref{pi}, we get
 \begin{align}\label{pii}
   \|\Pi(u,u)\|_{\widetilde{L}^1_t(H^{s_1})}\lesssim \|u\|_{\widetilde{L}^\infty_t(H^{s_1})}\|u\|_{\widetilde{L}^\infty_t(H^{s_1})}t,
 \end{align}
where we have used the fact that $H^{s_1}\hookrightarrow C^{0,1}.$\\
As for the term $\mathcal{P}\big((\nabla\cdot\xi^{m})\xi^{m}\big),$ by taking advantage of Bony's decomposition and of Lemmas \ref{T}-\ref{R}, we have
  \begin{align}
  \|\mathcal{P}\big((\nabla\cdot\xi^{m})\xi^{m}\big)\|_{H^{s_1}}
  \lesssim
  &\|\big((\nabla\cdot\xi^{m})\xi^{m}\big)\|_{H^{s_1}}    \\\nonumber
  \lesssim
  &\|\nabla\cdot\xi^{m}\|_{B^{s_1-(s_2+\frac{3}{2})}_{\infty,\infty}}
  \|\xi\|_{H^{s_2+\frac{3}{2}}}+\|\xi^{m}\|_{B^{s_1-(s_2+\frac{3}{2})}_{\infty,\infty}}
  \|\nabla\cdot\xi\|_{H^{s_2+\frac{3}{2}}}    \\\nonumber
  \lesssim&
  \|\nabla\cdot\xi^{m}\|_{H^{s_1-(s_2+\frac{3}{2})+\frac{d}{2}}}
  \|\xi\|_{H^{s_2+\frac{3}{2}}}+\|\xi^{m}\|_{H^{s_1-(s_2+\frac{3}{2})+\frac{d}{2}}}
  \|\nabla\cdot\xi\|_{H^{s_2+\frac{3}{2}}}     \\\nonumber
  \lesssim&
   \|\nabla\cdot\xi^{m}\|_{H^{s_2+\frac{1}{2}-\varepsilon}}
  \|\xi\|_{H^{s_2+\frac{3}{2}}}+\|\xi^{m}\|_{H^{s_2+\frac{1}{2}-\varepsilon}}
  \|\nabla\cdot\xi\|_{H^{s_2+\frac{3}{2}}}    \\ \nonumber \lesssim&\|\xi^{m}\|_{H^{s_2+\frac{3}{2}-\varepsilon_0}}
  \|\xi\|_{H^{s_2+\frac{5}{2}}},
  \end{align}
  where we have used $s_1-(s_2+\frac{3}{2})+\frac{d}{2}\leq \frac{d}{2}\leq s_1-1=s_2+\frac{1}{2}-\varepsilon,$ and $0<\varepsilon_0<\varepsilon.$\\
Inserting this inequality and (\ref{pii}) into (\ref{u'}), we get
\begin{align}\label{u}
  \|u^{m+1}\|_{\widetilde{L}^\infty_t(H^{s_1})}\lesssim & exp(C\int_0^t\|u^{m}\|_{H^{s_1}}dt')\Big(
  \|u_0\|_{H^{s_1}}+\|u^{m}\|_{\widetilde{L}^\infty_t(H^{s_1})}^2t
  \\\nonumber&+\|\xi^m\|_{\widetilde{L}^{\frac{4}{1-2\varepsilon_0}}(H^{s_2+1+\frac{2}{\frac{4}{1-2\varepsilon_0}}})}
    \|\xi^m\|_{\widetilde{L}^{\frac{4}{3}}(H^{s_2+1+\frac{2}{\frac{4}{3}}})}t^{\frac{\varepsilon_0}{2}}\Big).
\end{align}

As regards $z^{m+1},$ it follows from Lemma \ref{tds} that
\begin{align*}
  &\|z^{m+1}\|_{\widetilde{L}^\infty_t(H^{s_2})}
  +\|z^{m+1}\|_{\widetilde{L}^1_t(H^{s_2+2})}\\\lesssim
  &(1+t)\Big(\|z_0\|_{H^{s_2}}+\|\nabla\cdot(u^{m}z^{m})\|_{\widetilde{L}^1_t(H^{s_2})}
  +\|\nabla\cdot\big((\nabla\cdot\xi^{m})\xi^{m}\big)\|_{\widetilde{L}^1_t(H^{s_2})}\Big).
  \end{align*}
  According to Lemma \ref{s_2+1}, we get
  \begin{align}\label{1}
   \|\nabla\cdot(u^{m}z^{m})\|_{\widetilde{L}^1_t(H^{s_2})}
   \lesssim\|u^{m}z^{m}\|_{\widetilde{L}^1_t(H^{s_2+1})}&\lesssim
   \|u^{m}\|_{\widetilde{L}^\infty_t(H^{s_2+1})}
   \|z^{m}\|_{\widetilde{L}^2_t(H^{s_2+1})}t^{\frac{1}{2}}\\\nonumber
   &\lesssim\|u^{m}\|_{\widetilde{L}^\infty_t(H^{s_1})}
   \|z^{m}\|_{\widetilde{L}^2_t(H^{s_2+1})}t^{\frac{1}{2}},
 \end{align}
 \begin{align}\label{2}
   \|\nabla\cdot\big((\nabla\cdot\xi^{m})\xi^{m}\big)\|_{\widetilde{L}^1_t(H^{s_2})}
  \lesssim\|\big((\nabla\cdot\xi^{m})\xi^{m}\big)\|_{\widetilde{L}^1_t(H^{s_2+1})}
&\lesssim\|(\nabla\cdot\xi^{m})\|_{\widetilde{L}^2_t(H^{s_2+1})}
\|\xi^{m}\|_{\widetilde{L}^\infty_t(H^{s_2+1})}t^{\frac{1}{2}}\\\nonumber
&\lesssim
\|\xi^{m}\|_{\widetilde{L}^2_t(H^{s_2+2})}
\|\xi^{m}\|_{\widetilde{L}^\infty_t(H^{s_2+1})}t^{\frac{1}{2}}.
 \end{align}
 Thus, we conclude that
 \begin{align}\label{n}
  &\|z^{m+1}\|_{\widetilde{L}^\infty_t(H^{s_2})}
  +\|z^{m+1}\|_{\widetilde{L}^1_t(H^{s_2+2})}\\\nonumber\lesssim
  &(1+t)\Big(\|z_0\|_{H^{s_2}}+\|u^{m}\|_{\widetilde{L}^\infty_t(H^{s_1})}\| z^{m}\|_{\widetilde{L}^2_t(H^{s_2+1})}t^{\frac{1}{2}}
  +\|\xi^{m}\|_{\widetilde{L}^2_t(H^{s_2+2})}
\|\xi^{m}\|_{\widetilde{L}^\infty_t(H^{s_2+1})}t^{\frac{1}{2}}\Big).
  \end{align}

 Similarly, combining Lemma \ref{s_2+1} with Lemma \ref{tds} yields
  \begin{align}\label{p}
  &\|\xi^{m+1}\|_{\widetilde{L}^\infty_t(H^{s_2+1})}
  +\|\xi^{m+1}\|_{\widetilde{L}^1_t(H^{s_2+3})}\\\nonumber\lesssim
  &(1+t)\Big(\|\xi_0\|_{H^{s_2}}+\|\mathcal{L}\big(u^{m}(\nabla \cdot \xi^{m})\big)\|_{\widetilde{L}^1_t(H^{s_2+1})}
  +\|\mathcal{L}\big(z^m\xi^{m}\big)\|_{\widetilde{L}^1_t(H^{s_2+1})}\Big)\\\nonumber\lesssim
  &(1+t)\Big(\|\xi_0\|_{H^{s_2}}+\|u^{m}\|_{\widetilde{L}^\infty_t(H^{s_2+1})}\|\nabla \cdot \xi^{m}\|_{\widetilde{L}^2_t(H^{s_2+1})}t^{\frac{1}{2}}
  +\|z^m\|_{\widetilde{L}^2_t(H^{s_2+1})}\|\xi^{m}\|_{\widetilde{L}^\infty_t(H^{s_2+1})}t^{\frac{1}{2}}\Big)\\
  \nonumber\lesssim
  &(1+t)\Big(\|\xi_0\|_{H^{s_2}}+\|u^{m}\|_{\widetilde{L}^\infty_t(H^{s_1})}\| \xi^{m}\|_{\widetilde{L}^2_t(H^{s_2+2})}t^{\frac{1}{2}}
  +\|z^m\|_{\widetilde{L}^2_t(H^{s_2+1})}\|\xi^{m}\|_{\widetilde{L}^\infty_t(H^{s_2+1})}t^{\frac{1}{2}}\Big).
  \end{align}

  Denote $$E^{m}(t)\triangleq \|u^{m}\|_{\widetilde{L}^\infty_t(H^{s_1})}+\|z^{m}\|_{\widetilde{L}^\infty_t(H^{s_2})\cap \widetilde{L}^1_t(H^{s_2+2})}+
  \|\xi^{m}\|_{\widetilde{L}^\infty_t(H^{s_2+1})\cap \widetilde{L}^1_t(H^{s_2+3})},$$
  and
  $$E^0\triangleq\|u_0\|_{H^{s_1}}+\|z_0\|_{H^{s_2}}+
  \|\xi_0\|_{H^{s_2+1}}.$$
  By using interpolation and plugging the inequalities (\ref{n}) and (\ref{p}) into (\ref{u}) yield
  \begin{align*}
    E^{m+1}(t)\leq C\big(e^{CE^{m}(t)t}+1+t\big)\Big(E^0+\big(E^{m}(t)\big)^2\big(t+t^{\frac{1}{2}}+t^{\frac{\varepsilon_0}{2}}\big)\Big).
  \end{align*}
  Let us choose a positive $T_0\leq1$ such that $exp{(8C^2E^0T_0)}\leq 2$ and $T_0^{\frac{\varepsilon_0}{2}}\leq \frac{1}{192C^2E_0}.$ The induction hypothesis then implies that $$E^{m}(T_0)\leq 8CE^0.$$
\subsubsection{Second step: Convergence of the sequence}
  Let us fix some positive $T$ such that $T\leq T_0,$ and $(2CE^0)^4T\leq 1.$ We frist consider the case $s_1\neq 2+\frac{d}{2}.$

  By taking the difference between the equations for $u^{m+1}$ and $u^{m},$ one finds that
   \begin{align}\label{uu}
     &(u^{m+1}-u^m)_t+u^{m}\cdot \nabla (u^{m+1}-u^m)\\=\nonumber&(u^{m-1}-u^{m})\nabla u^{m}-\Pi(u^{m}-u^{m-1},u^{m}+u^{m-1})\\\nonumber&+\mathcal{P}\big((\nabla\cdot\xi^m)(\xi^{m}-\xi^{m-1})\big)
   +\mathcal{P}\big((\nabla\cdot\xi^m-\nabla\cdot\xi^{m-1})\xi^{m-1}\big).
   \end{align}
   Thanks to Lemma (\ref{deng}), we have
   \begin{align}\label{i}
     &\|\Pi(u^{m}-u^{m-1},u^{m}+u^{m-1})\|_{\widetilde{L}^1_t(H^{s_1-1})}
     \\\nonumber\lesssim &\|u^{m-1}-u^{m}\|_{\widetilde{L}^\infty_t(H^{s_1-1})}
     (\|u^{m}\|_{\widetilde{L}^\infty_t(H^{s_1})}
     +\|u^{m-1}\|_{\widetilde{L}^\infty_t(H^{s_1})})t.
   \end{align}
   From Lemmas \ref{T}-\ref{R}, we deduce that
   \begin{align}
   \|(u^{m-1}-u^{m})\nabla u^{m+1}\|_{\widetilde{L}^1_t(H^{s_1-1})}\lesssim \|u^{m-1}-u^{m}\|_{\widetilde{L}^\infty_t(H^{s_1-1})}
     \|u^{m-1}\|_{\widetilde{L}^\infty_t(H^{s_1})}t,\end{align}
     \begin{align}
  &\|\mathcal{P}\big((\nabla\cdot\xi^{m})(\xi^{m}-\xi^{m-1})\big)\|_{H^{s_1-1}}\\\nonumber
  \lesssim
  &\|\big((\nabla\cdot\xi^{m})(\xi^{m}-\xi^{m-1})\big)\|_{H^{s_1-1}}    \\\nonumber
  \lesssim
  &\|\nabla\cdot\xi^{m}\|_{B^{s_1-1-(s_2+\frac{3}{2})}_{\infty,\infty}}
  \|\xi^{m}-\xi^{m-1}\|_{H^{s_2+\frac{3}{2}}}+\|\xi^{m}-\xi^{m-1}\|_{B^{s_1-(s_2+\frac{3}{2})}_{\infty,\infty}}
  \|\nabla\cdot\xi^{m}\|_{H^{s_2+\frac{1}{2}}}    \\\nonumber
  \lesssim&
  \|\nabla\cdot\xi^{m}\|_{H^{s_1-1-(s_2+\frac{3}{2})+\frac{2}{d}}}
  \|\xi^{m}-\xi^{m-1}\|_{H^{s_2+\frac{3}{2}}}+\|\xi^{m}-\xi^{m-1}\|_{H^{s_1-(s_2+\frac{3}{2})+\frac{2}{d}}}
  \|\nabla\cdot\xi^m\|_{H^{s_2+\frac{1}{2}}}     \\\nonumber
  \lesssim&
   \|\nabla\cdot\xi^{m}\|_{H^{s_2-\frac{1}{2}-\varepsilon_0}}
  \|\xi^{m}-\xi^{m-1}\|_{H^{s_2+\frac{3}{2}}}+\|\xi^{m}-\xi^{m-1}\|_{H^{s_2+\frac{1}{2}-\varepsilon_0}}
  \|\nabla\cdot\xi^m\|_{H^{s_2+\frac{1}{2}}},
  \end{align}
  where we have used $s_1-1-(s_2+\frac{3}{2})+\frac{2}{d}\leq \frac{2}{d}-1\leq s_1-2\leq s_2+\frac{1}{2}-\varepsilon.$ Similarly,
     \begin{align}\label{iii}
  &\|\mathcal{P}\big((\nabla\cdot\xi^m-\nabla\cdot\xi^{m-1})\xi^{m-1}\big)\|_{H^{s_1-1}}
  \\\nonumber\lesssim&
   \|\nabla\cdot\xi^m-\nabla\cdot\xi^{m-1}\|_{H^{s_2-\frac{1}{2}-\varepsilon_0}}
  \|\xi^{m-1}\|_{H^{s_2+\frac{3}{2}}}+\|\xi^{m-1}\|_{H^{s_2+\frac{1}{2}-\varepsilon_0}}
  \|\nabla\cdot\xi^m-\nabla\cdot\xi^{m-1}\|_{H^{s_2+\frac{1}{2}}}.
  \end{align}
   Applying Lemma \ref{ts} to (\ref{uu}) thus yields
   \begin{align}\label{uug}
  \|u^{m+1}-&u^{m}\|_{\widetilde{L}^\infty_t(H^{s_1-1})}\lesssim exp(C\int_0^t\|u^{m}\|_{H^{s_1}}dt')\\\nonumber&\times\Big(\|u^{m-1}-u^{m}\|_{\widetilde{L}^\infty_t(H^{s_1-1})}
  (\|u^{m}\|_{\widetilde{L}^\infty_t(H^{s_1})}
     +\|u^{m-1}\|_{\widetilde{L}^\infty_t(H^{s_1})})t\\\nonumber&~~~+
  (\|\xi^m\|_{\widetilde{L}^\infty_t(H^{s_2+1})}+\|\xi^{m-1}\|_{\widetilde{L}^\infty_t(H^{s_2+1})})
    \|\xi^m-\xi^{m-1}\|_{\widetilde{L}^{\frac{4}{3}}_t(H^{s_2+\frac{2}{\frac{4}{3}}})}t^{\frac{1}{4}}
    \\\nonumber&~~~+\|\xi^m-\xi^{m-1}\|_{\widetilde{L}^{\frac{4}{1-2\varepsilon_0}}_t(B^{s_2+\frac{2}{\frac{4}{1-2\varepsilon_0}}})}
    (\|\xi^m\|_{\widetilde{L}^4_t(H^{s_2+1+\frac{2}{4}})}+\|\xi^{m-1}\|_{\widetilde{L}^4_t(H^{s_2+1+\frac{2}{4}})}t^{\frac{1+\varepsilon_0}{2}}\Big).
\end{align}

Note that
\begin{align*}
  (z^{m+1}-z^{m})_t-&\triangle \Big(z^{m+1}-z^{m})=-\nabla\cdot\big(u^m (z^{m}-z^{m-1})-(u^{m}-u^{m-1}) z^{m-1}\big)\\&+\nabla\big((\nabla\cdot\xi^m) (\xi^{m}-\xi^{m-1})\big)-\nabla\big((\nabla\cdot \xi^{m}-\nabla\cdot \xi^{m-1})\xi^{m-1} \big).
\end{align*}
By virtue of Lemma \ref{s_2+1}, we get
\begin{align}\label{44}
   \|\nabla\cdot\big(u^m (z^{m}-z^{m-1})\big)
    \|_{H^{s_2-1}}
   \lesssim&
   \|u^m (z^{m}-z^{m-1})
    \|_{H^{s_2}}\\\nonumber\lesssim&
   \|u^m\|_{H^{s_2+\frac{1}{2}}}\|z^{m}-z^{m-1}
    \|_{H^{s_2}}\lesssim
   \|u^m\|_{H^{s_1}}\|z^{m}-z^{m-1}
    \|_{H^{s_2}},
   \end{align}
   \begin{align}\label{55}
   \|\nabla\cdot\big(u^{m}-u^{m-1})z^{m-1}\big)
    \|_{H^{s_2-1}}\lesssim&\|(u^{m}-u^{m-1})z^{m-1}
    \|_{H^{s_2}}\\\nonumber\lesssim&
   \|u^{m}-u^{m-1}\|_{H^{s_2}}
   \|z^{m}\|_{H^{s_2+\frac{1}{2}}}\lesssim
   \|u^{m}-u^{m-1}\|_{H^{s_1-1}}
   \|z^{m}\|_{H^{s_2+\frac{1}{2}}},
   \end{align}
   \begin{align}\label{66}
   \|\nabla\big((\nabla\cdot\xi^m) (\xi^{m}-\xi^{m-1})\big)
    \|_{H^{s_2-1}}
    \lesssim&
   \|(\nabla\cdot\xi^m) (\xi^{m}-\xi^{m-1})\|_{H^{s_2}}\\\nonumber\lesssim&
   \|\nabla\cdot\xi^m\|_{H^{s_2+\frac{1}{2}}}
   \|\xi^{m}-\xi^{m-1}\|_{H^{s_2}},
   \end{align}
   \begin{align}\label{77}
   \|\nabla\big((\nabla\cdot \xi^{m}-\nabla\cdot \xi^{m-1})\xi^{m-1} \big)
    \|_{H^{s_2-1}}
  \lesssim&
   \|(\nabla\cdot \xi^{m}-\nabla\cdot \xi^{m-1})\xi^{m-1}\|_{H^{s_2}}\\\nonumber\lesssim&
   \|\nabla\cdot \xi^{m}-\nabla\cdot \xi^{m-1}\|_{H^{s_2+\frac{1}{2}}}
   \|\xi^{m-1}\|_{H^{s_2}}.
 \end{align}
 Hence Lemma \ref{tds} implies that
 \begin{align}\label{nng}
  &\|z^{m+1}-z^{m}\|_{\widetilde{L}^\infty_t(H^{s_2-1})}
  +\|z^{m+1}-z^{m}\|_{\widetilde{L}^1_t(H^{s_2+1})}\\\nonumber\lesssim
  &(1+t)\Big(\|u^{m}\|_{\widetilde{L}^\infty_t(H^{s_1})}
   \|z^{m}-z^{m-1}\|_{\widetilde{L}^2_t(H^{s_2})}t^{\frac{1}{2}}
   +\|u^{m}-u^{m-1}\|_{\widetilde{L}^\infty_t(H^{s_1-1})}
   \|z^{m}\|_{\widetilde{L}^4_t(H^{s_2+\frac{2}{4}}_{p_2,r_2})}t^{\frac{3}{4}}
  \\\nonumber&+\|\xi^{m}\|_{\widetilde{L}^4_t(H^{s_2+1+\frac{2}{4}})}
   \|\xi^{m}-\xi^{m-1}\|_{\widetilde{L}^\infty_t(H^{s_2})}t^{\frac{3}{4}}+
   \|\xi^{m}-\xi^{m-1}\|_{\widetilde{L}^{\frac{4}{3}}_t(H^{s_2+\frac{2}{\frac{4}{3}}}_{p_2,r_2})}
   \|\xi^{m-1}\|_{\widetilde{L}^\infty_t(H^{s_2+1})}t^{\frac{1}{4}}\Big).
  \end{align}
  Similarly, we get
  \begin{align}\label{ppg}
  &\|\xi^{m+1}-\xi^{m}\|_{\widetilde{L}^\infty_t(H^{s_2})}
  +\|\xi^{m+1}-\xi^{m}\|_{\widetilde{L}^1_t(H^{s_2+2})}\\\nonumber\lesssim
  &(1+t)\Big(\|\mathcal{L}\big(u^{m}\nabla\cdot(\xi^{m}-\xi^{m-1})\big)\|_{\widetilde{L}^1_t(H^{s_2})}
   +\|\mathcal{L}\big((u^{m}-u^{m-1})\nabla\cdot\xi^{m-1}\big)\|_{\widetilde{L}^1_t(H^{s_2})}\\\nonumber&
   +\|\mathcal{L}\big(z^{m}(\xi^{m}-\xi^{m-1})\big)\|_{\widetilde{L}^1_t(H^{s_2})}
  +\|\mathcal{L}\big((z^m-z^{m-1})\xi^{m}\big)\|_{\widetilde{L}^{1}_t(H^{s_2})}\Big)\\\nonumber\lesssim
  &(1+t)\Big(\|u^{m}\|_{\widetilde{L}^\infty_t(H^{s_2+\frac{1}{2}})}
   \|\nabla\cdot(\xi^{m}-\xi^{m-1})\|_{\widetilde{L}^2_t(H^{s_2})}t^{\frac{1}{2}}
   +\|u^{m}-u^{m-1}\|_{\widetilde{L}^\infty_t(H^{s_2})}
   \|\nabla\cdot\xi^{m-1}\|_{\widetilde{L}^4_t(H^{s_2+\frac{2}{4}})}t^{\frac{3}{4}}\\\nonumber&
   +\|z^{m}\|_{\widetilde{L}^\infty_t(H^{s_2})}\|\xi^{m}-\xi^{m-1}\|_{\widetilde{L}^4_t(H^{s_2+\frac{2}{4}})}t^{\frac{3}{4}}
  +\|z^m-z^{m-1}\|_{\widetilde{L}^2_t(H^{s_2})}
   \|\xi^{m}\|_{\widetilde{L}^{\infty}_t(H^{s_2+\frac{1}{2}})}t^{\frac{1}{2}}\Big)\\\nonumber\lesssim
  &(1+t)\Big(\|u^{m}\|_{\widetilde{L}^\infty_t(H^{s_1})}
   \|\xi^{m}-\xi^{m-1}\|_{\widetilde{L}^2_t(H^{s_2+1})}t^{\frac{1}{2}}
   +\|u^{m}-u^{m-1}\|_{\widetilde{L}^\infty_t(H^{s_1-1})}
   \|\xi^{m-1}\|_{\widetilde{L}^4_t(H^{s_2+1+\frac{2}{4}})}t^{\frac{3}{4}}\\\nonumber&
   +\|z^{m}\|_{\widetilde{L}^\infty_t(H^{s_2})}\|\xi^{m}-\xi^{m-1}\|_{\widetilde{L}^4_t(H^{s_2+\frac{2}{4}})}t^{\frac{3}{4}}
  +\|z^m-z^{m-1}\|_{\widetilde{L}^2_t(H^{s_2})}
   \|\xi^{m}\|_{\widetilde{L}^{\infty}_t(H^{s_2+1})}t^{\frac{1}{2}}\Big).
  \end{align}
   Denote \begin{align*}F^{m}(t)\triangleq \|u^{m+1}-u^{m}\|_{\widetilde{L}^\infty_t(H^{s_1-1})}&+\|z^{m+1}-z^{m}\|_{\widetilde{L}^\infty_t
   (H^{s_2-1})\cap\widetilde{L}^1_t(H^{s_2+1})}\\&+
  \|\xi^{m+1}-\xi^{m}\|_{\widetilde{L}^\infty_t(H^{s_2})
  \cap\widetilde{L}^1_t(H^{s_2+2})}.\end{align*}
  Plugging the inequalities (\ref{nng}) and (\ref{ppg}) into (\ref{uug}) yields
  \begin{align*}
    F^{m+1}(T)\leq & C\big(e^{CE^{m}(T)T}+1+T\big)\big(E^{m}(T)+E^{m-1}(T)\big)F^{m}(T)\big(T+T^{\frac{1+\varepsilon_0}{2}}+T^{\frac{1}{2}}+T^{\frac{3}{4}}+T^{\frac{1}{4}}\big)
    \\\leq &CE^0T^{\frac{1}{4}}F^{m}(T)\leq \frac{1}{2}F^{m}(T).
  \end{align*}
 Hence, $(u^m,z^m,\xi^m)|_{m\in \mathbb{N}}$ is a Cauchy sequence in  $\widetilde{L}^\infty_T(H^{s_1-1})\times \big(\widetilde{L}^\infty_T(H^{s_2-1})\cap\widetilde{L}^1_T(H^{s_2+1})\big)\times \big(\widetilde{L}^\infty_T(H^{s_2})\cap\widetilde{L}^1_T(H^{s_2+2})\big).$\\
 In the case $s_1=2+\frac{d}{2},$ for every $\zeta\in(0,1),$ we have
 $$s_1-\zeta>1+\frac{d}{2}, \textit{and}~  s_2-\zeta+\frac{3}{2}> s_1-\zeta\geq s_2-\zeta+1.$$
 Following along the same lines as above, we have
  $(u^m,z^m,\xi^m)|_{m\in \mathbb{N}}$ is a Cauchy sequence in  $\widetilde{L}^\infty_T(H^{s_1-\zeta-1})\times \big(\widetilde{L}^\infty_T(H^{s_2-\zeta-1})\cap\widetilde{L}^1_T(H^{s_2-\zeta+1})\big)\times \big(\widetilde{L}^\infty_T(H^{s_2-\zeta})\cap\widetilde{L}^1_T(H^{s_2-\varsigma+2})\big).$
  \subsubsection{Third step: Passing to the limit}
  Since the case $s_1=2+\frac{d}{2}$ works the same way, we only consider the case $s_1\neq2+\frac{d}{2}.$ Let $(u,z,\xi)$ be the limit of the sequence $(u^m,z^m,\xi^m)|_{m\in \mathbb{N}}.$ We see that $(u,z,\xi)\in\widetilde{L}^\infty_T(H^{s_1-1})\times \big(\widetilde{L}^\infty_T(H^{s_2-1})\cap\widetilde{L}^1_T(H^{s_2+1})\big)\times \big(\widetilde{L}^\infty_T(H^{s_2})\cap\widetilde{L}^1_T(H^{s_2+2})\big).$ Using Lemma \ref{Fadou} with the uniform bounds given in Step 1, we see that $(u,n,p)\in\widetilde{L}^\infty_T(H^{s_1})\times \big(\widetilde{L}^\infty_T(H^{s_2})\big)\times \big(\widetilde{L}^\infty_T(H^{s_2+1})\big).$ Next, by interpolating we discover that $(u^m,z^m,\xi^m)$ tends to $(u,z,\xi)$ in every space $\widetilde{L}^\infty_T(H^{s_1-\eta})\times \big(\widetilde{L}^\infty_T(H^{s_2-\eta})\cap\widetilde{L}^1_T(H^{s_2+1})\big)\times
  \big(\widetilde{L}^\infty_T(H^{s_2+1-\eta})\cap\widetilde{L}^1_T(H^{s_2+2})\big),$ with $\eta>0,$ which suffices to pass to the limit in the system (\ref{s2}).

  We still have to prove that $(z,\xi)\in \widetilde{L}^1_T(H^{s_2+2})\times\widetilde{L}^1_T(H^{s_2+3}).$ In fact, it is easy to check  that $(\partial_t z-\triangle z,\partial_t \xi-\triangle \xi)\in\widetilde{L}^1_T(H^{s_2})\times\widetilde{L}^1_T(H^{s_2+1}).$  Hence according to Lemma \ref{tds}, $(z,\xi)\in \widetilde{L}^1_T(H^{s_2+2})\times\widetilde{L}^1_T(H^{s_2+3}).$\\
  Finally, following along the same lines as in Theorem 3.19 of \cite{keben}, we can show that
   \begin{align}\label{lian}
    (u,z,\xi)\in C([0,T];H^{s_1})\times C([0,T];H^{s_2})\times C([0,T];H^{s_2+1}).
  \end{align}
\subsection{Uniqueness for the system (\ref{s2})}
 Without loss of generality, we may assume that $s_1<2+\frac{d}{2}.$ Assume that we are given $(u_1,z_1,\xi_1)$ and $(u_2,z_2,\xi_2)$, two solutions of the system (\ref{s2}) (with the same initial data) satisfying the regularity assumptions of Theorem \ref{a1}. In order to show these two solutions coincide,  we first denote
$$E^i(T)\triangleq
  \|u_i\|_{\widetilde{L}^\infty_T(H^{s_1})}+\|z_i\|_{\widetilde{L}^\infty_T(H^{s_2})}+
  \|\xi_i\|_{\widetilde{L}^\infty_t(H^{s_2+1})},~~i=1,2,$$
$$F(t_1,t_2)\triangleq \| u_2-u_1\|_{\widetilde{L}^\infty_{[t_1,t_2]}(H^{s_1-1})}+\|z_2-z_1 \|_{\widetilde{L}^\infty_{[t_1,t_2]}(H^{s_2-1})}+
  \|\xi_2-\xi_1\|_{\widetilde{L}^\infty_{[t_1,t_2]}(H^{s_2}) },$$
and  $$T_0=\sup\{0\leq T'\leq T~|~(u_1,z_1,\xi_1)=
(u_2,z_2,\xi_2)~on ~[0,T']\}. $$
 We deduce from the definition of $T_0$ and the continuity of $(u_i,z_i,\xi_i)$ that $$\big(u_1(T_0),z_1(T_0),\xi_1(T_0)\big)=
\big(u_2(T_0),z_2(T_0),\xi_2(T_0)\big).$$

If $T_0<T,$ repeating the same arguments as we were used for the proof of the convergence of the approximate solutions in the above subsection, we get
\begin{align*}
    F(T_0+\widetilde{T})\leq  C\big(e^{CE^{2}(T)\widetilde{T}}+1+\widetilde{T}\big)\big(E^{2}(T)
    +E^{1}(T)\big)F(T_0+\widetilde{T})\big(\widetilde{T}
    +\widetilde{T}^{\frac{1+\varepsilon_0}{2}}
    +\widetilde{T}^{\frac{1}{2}}+\widetilde{T}^{\frac{3}{4}}+\widetilde{T}^{\frac{1}{4}}\big)
    .
  \end{align*}
  We conclude that $F(T_0+\widetilde{T})=0$ with sufficiently small $\widetilde{T}.$ Thus, $(u_1,z_1,\xi_1)=
(u_2,z_2,\xi_2)$ on $[T_0,T_0+\widetilde{T}]$, which stands in contradiction to the definition of $T_0.$ Hence $T_0=T,$
  and the proof of uniqueness is completed.
\subsection{Properties of $(u,z,\xi)$}
\subsubsection{$\nabla\cdot u=0$}
  Suppose that $(u,z,\xi)$ satisfies the system (\ref{s2}) in $\in \widetilde{L}^\infty_{T}(H^{s_1}(\mathbb{R}^d))\times \big(\widetilde{L}^\infty_{T}(H^{s_2}(\mathbb{R}^d))
\cap\widetilde{L}^1_{T}(H^{s_2+2}(\mathbb{R}^d))\big)\times \big(\widetilde{L}^\infty_{T}(H^{s_2+1}(\mathbb{R}^d))
\cap\widetilde{L}^1_{T}(H^{s_2+3}(\mathbb{R}^d))\big).$  We check that $u$ is divergence free. This may be achieved by applying $\nabla\cdot$ to the first equation of the system (\ref{s2}). Denote $s_1'=s_1-1$, if $s_1\neq 2+\frac{d}{2};$ and $s_1'=s_1-\zeta-1$, for some $\zeta\in(0,1),$ if $s_1=2+\frac{d}{2}.$  We get
  \begin{align*}
    (\partial_t+u\cdot \nabla)(\nabla\cdot u)=-\nabla\cdot\Pi(u,u)-tr(Du)^2.
  \end{align*}
  Lemma \ref{ts} and Lemma \ref{yyy} ensure that
  \begin{align*}
    \|\nabla\cdot u\|_{H^{s_1'}}\lesssim&
    \int_0^texp\Big(C\int_{t'}^t\|u\|_{H^{s_1}}dt''\Big)\|\nabla\cdot\Pi(u,u)+tr(Du)^2\|_{H^{s'_1}}dt'\\ \lesssim&
    \int_0^texp\Big(C\int_{t'}^t\|u\|_{H^{s_1}}dt''\Big)\|\nabla\cdot u\|_{B^{0}_{\infty,\infty}}\|u\|_{H^{s'_1+1}}dt',\\
    \lesssim&
    \int_0^texp\Big(C\int_{t'}^t\|u\|_{H^{s_1}}dt''\Big)\|\nabla\cdot u\|_{H^{s_1'}}\|u\|_{H^{s_1}}dt',
  \end{align*}
  where we have used $H^{s_1'}\hookrightarrow B^{0}_{\infty,\infty},$ and $H^{s_1'+1} \hookrightarrow H^{s_1}.$ Using Gronwall's inequality, we conclude that $\nabla\cdot u=0.$
 \subsubsection{$\mathcal{L}\xi=\xi$}
  Note that
\begin{align}\label{xixi}
 \xi = e^{t\triangle}\xi_0-\int_0^te^{(t-t')\triangle}\mathcal{L}\big(u(\nabla\cdot \xi)+z\xi\big)dt'.
\end{align}
Applying $\mathcal{L}$ to the above equation yields
\begin{align*}
\mathcal{L} \xi = e^{t\triangle}\mathcal{L}\xi_0-\int_0^te^{(t-t')\triangle}\mathcal{L}^2\big(u(\nabla\cdot \xi)+z\xi\big)dt'.
\end{align*}
It is easy to check that
\begin{align*}
  \mathcal{L}\xi_0=-\nabla (-\triangle)^{-1}\nabla\cdot(-\nabla (-\triangle)^{-1}a_0)=\nabla (-\triangle)^{-1}(\nabla\cdot\nabla ) (-\triangle)^{-1}a_0=-\nabla (-\triangle)^{-1}a_0=\xi_0,\\
  \mathcal{L}^2=-\nabla (-\triangle)^{-1}\nabla\cdot(-\nabla (-\triangle)^{-1}\nabla\cdot)=\nabla (-\triangle)^{-1}(\nabla\cdot\nabla)(-\triangle)^{-1}\nabla\cdot=-\nabla (-\triangle)^{-1}\nabla\cdot=\mathcal{L}.
\end{align*}
Hence, $\mathcal{L}\xi=\xi.$
\subsubsection{Nonnegative of $z\pm\nabla\cdot \xi$}
Let $a=\frac{z+\nabla \cdot \xi}{2},$ and $b=\frac{z-\nabla \cdot \xi}{2}.$ As $\nabla \cdot\mathcal{L}=\nabla \cdot$ and $\nabla\cdot u=0,$ one finds that $(a,b)$ solves the following system:
\begin{align}\tag{ab}\label{ab}
\left\{
\begin{array}{l}
a_t+u\cdot \nabla a-\triangle a=-\nabla\cdot(a\xi),  \\[1ex]
b_t+u\cdot \nabla b- \triangle b=\nabla\cdot(b\xi), \\[1ex]
\nabla \cdot \xi=a-b,\\[1ex]
(a,b)|_{t=0}=(\frac{z_0+\nabla \cdot \xi_0}{2},\frac{z_0-\nabla \cdot \xi_0}{2}),
\end{array}
\right.
\end{align}
  We test the first equation of the system (\ref{ab}) with $(a^-)\triangleq\sup \{-a,0\}.$ After integrating by parts, we obtain
  \begin{align*}
    \frac{1}{2}\frac{d}{dt}\|a^-\|_{L^{2}}^{2}+\|\nabla a^-\|_{L^{2}}^{2}=-\int_{\mathbb{R}^d}\frac{1}{2}(\nabla \cdot \xi)(a^-)^2dx\leq \frac{1}{2}\|\nabla \cdot \xi\|_{L^\infty}\|a^-\|_{L^2}^2.
  \end{align*}
  Gronwall's Lemma implies that
  \begin{align*}
  \|a^-\|_{L^{2}}^{2}\leq \|a^-_0\|_{L^{2}}^{2}exp\{\|\nabla \cdot \xi\|_{L^1_t(L^\infty)}\}.
  \end{align*}
 Since $a_0\geq 0,$ and $$\|\nabla \cdot \xi\|_{L^1_t(L^\infty)}\lesssim
 \|\nabla \cdot \xi\|_{L^1_t(H^{s_2+\frac{1}{2}})}\lesssim
 \|\nabla \cdot \xi\|_{L^2_t(H^{s_2+1})}t^{\frac{1}{2}},$$
 we have $\|a^-\|_{L^{2}}^{2}=0.$ Hence $a
\geq 0,$ a.e. on $[0,T]\times \mathbb{R}^d.$
  Repeating the same steps for $b$ implies $b\geq 0,~a.e.~on ~[0,T]\times \mathbb{R}^d,$
  and thus $z=a+b\geq 0,~a.e.~on ~[0,T]\times \mathbb{R}^d.$

\subsection{A global existence result in dimension $d=2$ }
According to the above subsections, local existence in $\widetilde{L}^\infty_{T}(H^{s_1}(\mathbb{R}^d))\times \big(\widetilde{L}^\infty_{T}(H^{s_2}(\mathbb{R}^d))
\cap\widetilde{L}^1_{T}(H^{s_2+2}(\mathbb{R}^d))\big)\times \big(\widetilde{L}^\infty_{T}(H^{s_2+1}(\mathbb{R}^d))
\cap\widetilde{L}^1_{T}(H^{s_2+3}(\mathbb{R}^d))\big)$ has already been proven. So we denote by $T^*$ the maximal time of existence of $(u,z,\xi).$ Suppose that $T^*$ is finite, under the assumption of \ref{jibentiaojian}, and assume further that $d=2$ and $s_2>1$, we have the following lemmas.
\subsubsection{Some useful lemmas}
\begin{lemm}\label{1l}
$\forall t\in[0,T^*),$ we have
\begin{align}\label{l1}
\|u(t)\|_{L^2}^2+\|\xi(t)\|_{L^2}^2+\int_0^t\|\nabla \xi\|_{L^2}^2dt'\lesssim\|u_0\|_{L^2}^2+\|\xi_0\|_{L^2}^2.
\end{align}
\end{lemm}
\noindent{Proof.} Multiplying the first and the third equations of the system (\ref{s2}) by $u$ and $\xi$ respectively, and integrating over $\mathbb{R}^d:$
\begin{align*}
  &\frac{1}{2}\frac{d}{dt}\|u\|_{L^2}^2
  =\int_{\mathbb{R}^2}(\nabla\cdot\xi)\xi udx,\\
  &\frac{1}{2}\frac{d}{dt}\|\xi\|_{L^2}^2+\int_{\mathbb{R}^2}(\nabla\cdot\xi)\xi udx+\|\nabla \xi\|_{L^2}^2=\int_{\mathbb{R}^d}-z|\xi|^2dx,
\end{align*}
where we have used
\begin{align*}
&\int_{\mathbb{R}^2}(u\cdot\nabla u)udx=
\int_{\mathbb{R}^2}-\frac{1}{2}(\nabla \cdot u)|u|^2dx
=0,\\
&\int_{\mathbb{R}^2}\Pi(u,u)udx=
\int_{\mathbb{R}^2}\nabla P_{\Pi}(u,u)udx
=\int_{\mathbb{R}^2}- P_{\Pi}(u,u)(\nabla \cdot u)dx=0,\\
&\int_{\mathbb{R}^2}\mathcal{P}\big((\nabla\cdot\xi)\xi\big)udx=
\int_{\mathbb{R}^2}(\nabla\cdot\xi)\xi(\mathcal{P}u)dx
=\int_{\mathbb{R}^2}(\nabla\cdot\xi)\xi udx,\\
&\int_{\mathbb{R}^2}\mathcal{L}\big(u(\nabla\cdot\xi)\big)\xi dx=
\int_{\mathbb{R}^2}u(\nabla\cdot\xi)(\mathcal{L}\xi) dx=
\int_{\mathbb{R}^2}u(\nabla\cdot\xi)\xi dx,\\
&\int_{\mathbb{R}^2}\mathcal{L}(-z\xi)\xi=
\int_{\mathbb{R}^2}-z\xi(\mathcal{L}\xi)=
\int_{\mathbb{R}^2}-z|\xi|^2,
\end{align*}
with $P_{\pi}(u,u)$ defined as in Lemma \ref{pi}.
 Summing the above equations and using the fact $z\geq 0,$ we find
 \begin{align*}
  \frac{1}{2}\frac{d}{dt}(\|u\|_{L^2}^2+\|\xi\|_{L^2}^2)+\|\nabla \xi\|_{L^2}^2\leq 0,
\end{align*}
from which it follows that (\ref{l1}) holds. \qed
 \begin{lemm}\label{2l}
$\forall t\in[0,T^*),$ $2\leq q\leq \infty,$ we have
\begin{align*}
  &\|z\pm\nabla \cdot \xi(t)\|_{L^q}\lesssim \|z_0+\nabla \cdot \xi_0\|_{L^q}+\|z_0-\nabla \cdot \xi_0\|_{L^q},\\
&\|z\pm\nabla \cdot \xi(t)\|_{L^2}^2+\int_0^t\|\nabla (z\pm\nabla \cdot \xi)\|_{L^2}^2dt'\lesssim \|z_0+\nabla \cdot \xi_0\|_{L^2}^2+\|z_0-\nabla \cdot \xi_0\|_{L^2}^2.
\end{align*}
\end{lemm}
 \noindent{Proof.} Since $s_2>1,$ (\ref{lian}) implies that $(a,b)\in \big(C([0,T];H^{s_2})\big)^2\hookrightarrow \big(C([0,T];L^q)\big)^2,$ with $2\leq q\leq \infty.$ By multiplying both sides of the first equation of the (\ref{ab}) system by $| a|^{p-2} a$ with $2\leq p< \infty,$ and integrating over $[0,t]\times\mathbb{R}^d,$ we get
 \begin{align}\label{n1}
   \frac{1}{p}\|a(t)\|_{L^p}^p+(p-1)\int_{\mathbb{R}^2}|a|^{p-2}|\nabla a|^2dx\leq\frac{1}{p}\|a_0\|_{L^p}^p-\frac{p-1}{p}
   \int_0^t\int_{\mathbb{R}^2}\nabla\cdot \xi|a|^pdxdt',
 \end{align}
 where we have used the estimates
 \begin{align*}
\int_{\mathbb{R}^2}| a|^{p-2} a(u\cdot \nabla a)dx=&
-\frac{1}{p}\int_{\mathbb{R}^2}(\nabla\cdot u)|a|^pdx=0,\\
-\int_{\mathbb{R}^2}| a|^{p-2} a\nabla\cdot(a\xi)dx=&
-\int_{\mathbb{R}^2}\xi\frac{1}{p}\nabla|a|^pdx-\int_{\mathbb{R}^2}(\nabla\cdot\xi)|a|^pdx\\=&
-\int_{\mathbb{R}^2}\frac{p-1}{p}(\nabla\cdot \xi)|a|^pdx.
 \end{align*}
 Repeating the same steps for $b$ yields
\begin{align}\label{p1}
   \frac{1}{p}\|b(t)\|_{L^p}^p+(p-1)\int_{\mathbb{R}^2}|b|^{p-2}|\nabla b|^2dx\leq\frac{1}{p}\|b_0\|_{L^p}^p+\frac{p-1}{p}
   \int_0^t\int_{\mathbb{R}^2}(\nabla\cdot\xi)|b|^pdxdt'.
 \end{align}
Adding up $(\ref{n1})$ and (\ref{p1}), we get
\begin{align*}
 &\frac{1}{p}(\|a(t)\|_{L^p}^p+\|b(t)\|_{L^p}^p)+(p-1)\int_{\mathbb{R}^2}|a|^{p-2}|\nabla a|^2dx+(p-1)\int_{\mathbb{R}^2}|b|^{p-2}|\nabla b|^2dx\\\leq &\frac{1}{p}(\|a_0\|_{L^p}^p+\|b_0\|_{L^p}^p)+\frac{p-1}{p}
   \int_0^t\int_{\mathbb{R}^2}\nabla\cdot\xi(|b|^p-|a|^p)dxdt'\\
   \leq &\frac{1}{p}(\|a_0\|_{L^p}^p+\|b_0\|_{L^p}^p)+\frac{p-1}{p}
   \int_0^t\int_{\mathbb{R}^2}(a-b)(b^p-a^p)dxdt'\\
   \leq &\frac{1}{p}(\|a_0\|_{L^p}^p+\|b_0\|_{L^p}^p),
\end{align*}
where we have used the non-negativity of $a,b.$
This thus leads to
\begin{align*}
&\|a(t)\|_{L^p}+\|b(t)\|_{L^p}\leq 2( \|a_0\|_{L^p}+\|b_0\|_{L^p}),\\
&\|a(t)\|_{L^2}^2+\|b(t)\|_{L^2}^2+\int_0^t(\|\nabla a\|_{L^2}^2+\|\nabla b\|_{L^2}^2)dt'\leq\|a_0\|_{L^2}^2+\|b_0\|_{L^2}^2.
\end{align*}
Passing to the limit as $p$ tends to infinite gives
\begin{align*}
\|a(t)\|_{L^\infty}+\|b(t)\|_{L^\infty}\leq 2( \|n_0\|_{L^\infty}+\|b_0\|_{L^\infty}).
\end{align*}
This completes the proof of the lemma. \qed
\begin{lemm}\label{3l}
$\forall t\in[0,T^*),$ we have
\begin{align}\label{l3}
\|\xi(t)\|_{L^\infty}\leq C(T^*)<\infty.
\end{align}
\end{lemm}
{\noindent Proof.}
It is easy to obtain from \ref{xixi} that
\begin{align*}
  \|\xi(t)\|_{L^{\infty}}\lesssim&\|\mathcal{F}^{-1}(e^{-t|x|^2})\|_{L^1} \|\xi_0\|_{L^{\infty}}+ \int_0^t\|\mathcal{F}^{-1}(e^{(t'-t)|x|^2})\|_{L^{2}}\Big(\|u(t')\|_{L^{2}}\| \nabla\cdot\xi(t')\|_{L^{\infty}}\\&+\|z(t')\|_{L^{\infty}}\|\xi(t')\|_{L^{2}}\Big)dt'\\\lesssim&
  \|\xi_0\|_{L^{\infty}}+\int_0^t(t'-t)^{-\frac{1}{2}}\Big(\|u(t')\|_{L^{2}}\| \nabla\cdot\xi(t')\|_{L^{\infty}}dt'+\|z(t')\|_{L^{\infty}}\|\xi(t')\|_{L^{2}}\Big)\\\lesssim&
  \|\xi_0\|_{L^{\infty}}+t^{\frac{1}{2}}\Big(\|u\|_{L^\infty_t(L^{2})}\| \nabla\cdot\xi\|_{L^\infty_t(L^{\infty})}+\|z\|_{L^\infty_t(L^{\infty})}\|\xi\|_{L^\infty_t(L^{2})}\Big).
 \end{align*}
 Applying Lemma \ref{1l} and Lemma \ref{2l} completes the proof. \qed

 Denote $w=\partial_2u_1-\partial_1u_2.$ Note that $\Pi(u,u)=\nabla P_{\Pi}(u,u),$ where $P_{\Pi}(u,u)$ defined as in Lemma \ref{pi}, $(\mathcal{P}-Id)\big((\nabla\cdot\xi)\xi\big)=\nabla (-\triangle)^{-1}\nabla\cdot\big((\nabla\cdot\xi)\xi\big),$ and $\xi=\mathcal{L}\xi=\nabla (-\triangle)^{-1}\nabla\cdot \xi.$ Then $w$ satisfies
 \begin{align}\label{lw}
   w_t+u\cdot \nabla w=\partial_2(\nabla\cdot\xi)\xi_1-\partial_1(\nabla\cdot\xi)\xi_2.
 \end{align}
 \begin{lemm}\label{uw}\cite{keben}
For all $s\in \mathbb{R}$ and $1\leq p, r \leq \infty,$ there exists a constant $C$ such that
\begin{align}\label{wu}
\|(Id-\triangle_{-1})u\|_{B^s_{p,r}}
\leq  C\|w\|_{B^s_{p,r}}.
\end{align}
\end{lemm}
 \begin{lemm}\label{4l}
$\forall t\in[0,T^*),$ we have
\begin{align}\label{l4}
\|\nabla u(t)\|_{L^2}\leq C(T^*)<\infty.
\end{align}
\end{lemm}
{\noindent Proof.}
 Multiplying (\ref{uw}) by $w$ and integrating over $\mathbb{R}^2:$
\begin{align*}
    \frac{1}{2}\frac{d}{dt}\|w\|_{L^{2}}^{2}\lesssim\|\nabla(\nabla\cdot\xi)\|_{L^2}\|\xi\|_{L^\infty}\|w\|_{L^2}.
  \end{align*}
   The Gronwall lemma implies that
  \begin{align*}
   \|w(t)\|_{L^{2}}\lesssim& \|w_0\|_{L^{2}}+\int_0^t\|\nabla(\nabla\cdot\xi)\|_{L^2}\|\xi\|_{L^\infty}dt'
   \lesssim\|w_0\|_{L^{2}}+\|\nabla(\nabla\cdot\xi)\|_{L^2_t(L^2)}\|\xi\|_{L^{\infty}_t(L^\infty)}t^{\frac{1}{2}}.
  \end{align*}
  Applying Lemma \ref{2l} and Lemma \ref{3l}, we have
  \begin{align}\label{www}
    \|w(t)\|_{L^2}\leq C(T^*)<\infty.
  \end{align}
  Next by splitting $u$ into low and high frequencies and using Lemma \ref{uw}, we see that
  \begin{align*}
   \|\nabla u\|_{L^{2}}\lesssim \|\triangle_{-1}\nabla u\|_{L^{2}}+\|(Id-\triangle_{-1})\nabla u\|_{L^{2}}\lesssim\| u\|_{L^{2}}+\|w\|_{L^{2}}.
  \end{align*}
 Applying Lemma \ref{2l} and the inequality \ref{www} then completes the proof of the lemma.
\begin{lemm}\label{5l}
\begin{align}\label{l5}
 \int_0^{T^*}\|\nabla(z\pm\nabla \cdot \xi)\|_{L^\infty}dt'<\infty.
 \end{align}
 \end{lemm}
{\noindent Proof.} First combining Lemma \ref{1l} and Lemma \ref{4l} with the Sobolev imbedding theorem, we see that
\begin{align}\label{uq}
u\in L^\infty_{T^*}(H^1)\hookrightarrow L^\infty_{T^*}(L^p),
\end{align}
with $2\leq p<\infty.$
Then we denote from the system (\ref{ab}) that
\begin{align*}
 \nabla a = \nabla e^{t\triangle}a_0-\int_0^te^{(t-t')\triangle}\nabla\big(u\cdot\nabla a+\xi\cdot\nabla a+a(\nabla\cdot\xi)\big)dt'.
\end{align*} We have
\begin{align*}
  &\|(\nabla a)(\tau)\|_{L^{\infty}}\\\lesssim&\| \mathcal{F}^{-1}(e^{-\tau|x|^2}x)\|_{L^1}\|a_0\|_{L^{\infty}}+ \int_0^\tau\|\mathcal{F}^{-1}(e^{(t'-\tau)|x|^2}x)\big\|_{L^{{q'}}}\|u(t')\|_{L^{q}}\| \nabla a(t')\|_{L^{\infty}}dt'\\
  &+\int_0^\tau\big\|\mathcal{F}^{-1}(e^{(t'-\tau)|x|^2}x)\big\|_{L^1}\big(\|\xi\|_{L^\infty}\|\nabla a\|_{L^{\infty}}+\|\nabla\cdot\xi\|_{L^\infty}\|a\|_{L^{\infty}}\big)(t')dt'\\\lesssim& \tau^{-\frac{1}{2}}\|a_0\|_{L^{\infty}}+ \int_0^\tau\frac{1}{{(\tau-t')}^{\frac{1}{2}+\frac{1}{q}}}\|u(t')\|_{L^{q}}\| \nabla a(t')\|_{L^{\infty}}dt'\\
  &+\int_0^\tau\frac{1}{{(\tau-t')}^{\frac{1}{2}}}\big(\|\xi\|_{L^\infty}\|\nabla a\|_{L^{\infty}}+\|\nabla\cdot\xi\|_{L^\infty}\|a\|_{L^{\infty}}\big)(t')dt'\\
  \lesssim&\tau^{-\frac{1}{2}}\| a_0\|_{L^{\infty}}+\tau^{\frac{1}{2}}\|\nabla\cdot\xi\|_{L^\infty_T(L^\infty)}\|a\|_{L^\infty_T(L^{
  \infty})}\\&+\int_0^\tau\big(\frac{1}{{(\tau-t')}^{\frac{1}{2}+\frac{1}{q}}}
  \|u(t')\|_{L^{q}}+\frac{1}{{(\tau-t')}^{\frac{1}{2}}}\| \xi(t')\|_{L^\infty}\big)\|\nabla a(t')\|_{L^{\infty}}dt'\\
  \lesssim&\tau^{-\frac{1}{2}}\| a_0\|_{L^{\infty}}+\tau^{\frac{1}{2}}\|\nabla\cdot\xi\|_{L^\infty_T(L^\infty)}\|a\|_{L^\infty_T(L^{
  \infty})}\\&
  +\int_0^\tau\big(\delta_1\frac{1}{{(\tau-t')}^{(\frac{1}{2}+\frac{1}{q})\gamma_1}}
  +\delta_2\frac{1}{{(\tau-t')}^{\frac{1}{2}\gamma_2}}\big)\|\nabla a(t')\|_{L^{\infty}}dt'\\&+\int_0^\tau\big(C_{\delta_1}
  \|u(t')\|^{\gamma'_1}_{L^{q}}+C_{\delta_2}\|\xi\|^{\gamma'_2}_{L^\infty}\big)\|\nabla a(t')\|_{L^{\infty}}dt',
 % &\tau^{-\frac{1}{2}}\| n_0\|_{L^{a_1}}+\int_0^\tau\frac{1}{{(\tau-t')}^{\frac{1}{2}}}\Big((\|u\|_{L^\infty}+\|\nabla \phi\|_{L^\infty})\|\nabla n\|_{L^{a_1}}+\|\triangle \phi\|_{L^\infty}\|n\|_{L^{a_1}}\Big)(t')dt'\\
%   \lesssim&\tau^{-\frac{1}{2}}\| n_0\|_{L^{a_1}}+\tau^{\frac{1}{2}}\|\triangle \phi\|_{L^\infty_\tau(L^\infty)}\|n\|_{L^\infty_\tau(L^{a_1})}+\delta
%   \int_0^\tau\frac{1}{{(\tau-t')}^{\frac{\gamma}{2}}}\|n\|_{L^{a_1}}(t')dt'\\
%   &+C_{\delta}\int_0^\tau(\|u\|_{L^\infty}^{\gamma'}+\|\nabla \phi\|_{L^\infty}^{\gamma'})\|n\|_{L^{a_1}}(t')dt',
\end{align*}
 with $2<q<\infty,$ $\frac{1}{q}+\frac{1}{q'}=1,$
  $ (\frac{1}{2}+\frac{1}{q})\gamma_1<1,~~\frac{1}{2}\gamma_2<1.
$
By means of the Young inequality for the time integral, we obtain,
\begin{align*}
  \|\nabla a\|_{L^1_t(L^{\infty})}\lesssim &t^{\frac{1}{2}}\| a_0\|_{L^{\infty}}+t^{\frac{3}{2}}\|\nabla\cdot\xi\|_{L^\infty_T(L^\infty)}\|a\|_{L^\infty_T(L^{\infty})}\\&
  +\big(\delta_1\frac{1}{1-(\frac{1}{2}+\frac{1}{q})\gamma_1}t^{1-(\frac{1}{2}+\frac{1}{q})\gamma_1}
  +\delta_2\frac{1}{1-\frac{1}{2}\gamma_2}t^{1-\frac{1}{2}\gamma_2}\big)\|\nabla a\|_{L^1_t(L^{a_1})}\\&+\int_0^t\int_0^\tau\big(C_{\delta_1}
  \|u(t')\|^{\gamma'_1}_{L^{q}}+C_{\delta_2}\|\xi(t')\|^{\gamma'_2}_{L^\infty}\big)\|\nabla a(t')\|_{L^{\infty}}dt'd\tau.
\end{align*}
Choosing $\delta_1\frac{1}{1-(\frac{1}{2}+\frac{1}{q})\gamma_1}(T^*)^{1-(\frac{1}{2}+\frac{1}{q})\gamma_1}
  +\delta_2\frac{1}{1-\frac{1}{2}\gamma_2}(T^*)^{1-\frac{1}{2}\gamma_2}=c\frac{1}{2}$ yields
\begin{align*}
   \|\nabla a\|_{L^1_t(L^{\infty})}\lesssim &(T^*)^{\frac{1}{2}}\| a_0\|_{L^{\infty}}+(T^*)^{\frac{3}{2}}\|\nabla\cdot\xi\|_{L^\infty_{T^*}(L^\infty)}\|a\|_{L^\infty_{T^*}(L^{\infty})}\\&
  +\big(C_{\delta_1}
  \|u(t')\|^{\gamma'_1}_{L^\infty_{T^*}(L^{q})}+C_{\delta_2}\|\xi\|^{\gamma'_2}_{L^\infty_{T^*}(L^\infty)}\big)\int_0^t\int_0^\tau\|\nabla a(t')\|_{L^{\infty}}dt'd\tau.
   \end{align*}
Gronwall's lemma thus implies that
\begin{align*}
\|\nabla a\|_{L^1_t(L^{\infty})}\leq &C\Big((T^*)^{\frac{1}{2}}\| a_0\|_{L^{\infty}}+(T^*)^{\frac{3}{2}}\|\nabla\cdot\xi\|_{L^\infty_{T^*}(L^\infty)}\|a\|_{L^\infty_{T^*}(L^{\infty})}\Big)\\ &\times exp\Big(\big(C_{\delta_1}
  \|u(t')\|^{\gamma'_1}_{L^\infty_{T^*}(L^{q})}+C_{\delta_2}\| \xi\|^{\gamma'_2}_{L^\infty_{T^*}(L^\infty)}\big)t\Big).
  \end{align*}
  Hence, Lemma \ref{2l}, Lemma \ref{3l} and the inequality (\ref{uq}) imply that
  \begin{align*}
 \int_0^{T^*}\|\nabla a\|_{L^\infty}dt'<\infty.
 \end{align*}
 Similar arguments for $b$ yield
 \begin{align*}
 \int_0^{T^*}\|\nabla b\|_{L^\infty}dt'<\infty.
 \end{align*}
 Therefore, the inequality (\ref{l5}) holds true. \qed
 \begin{lemm}\label{6l}
$\forall t\in[0,T^*),$ we have
\begin{align}\label{l6}
\|u(t)\|_{B^1_{\infty,\infty}}\leq C(T^*)<\infty.
\end{align}
\end{lemm}
{\noindent Proof.} we deduce from the inequality (\ref{lw}) that
 \begin{align*}
   \|w(t)\|_{L^{\infty}}\lesssim& \|w_0\|_{L^{\infty}}+\int_0^t\|\nabla(\nabla\cdot\xi)\|_{L^\infty}\|\xi\|_{L^\infty}dt'
   \lesssim\|w_0\|_{L^{\infty}}+\|\nabla(\nabla\cdot\xi)\|_{L^1_t(L^\infty)}\|\xi\|_{L^{\infty}_t(L^\infty)}.
  \end{align*}
  By splitting $u$ into low and high frequencies and using Lemma \ref{uw}, we see that
  \begin{align*}
   \| u\|_{B^{1}_{\infty,\infty}}\lesssim& \|\triangle_{-1} u\|_{L^{\infty}}+\|(Id-\triangle_{-1}) \nabla u\|_{B^{0}_{\infty,\infty}}\\\lesssim&
  \|u\|_{L^{2}}+\|w\|_{B^{0}_{\infty,\infty}}\lesssim\|u\|_{L^{2}}+\|w\|_{L^\infty}.
  \end{align*}
Applying Lemma \ref{1l}, Lemma \ref{3l} and Lemma \ref{5l} completes the proof of the lemma. \qed
  \subsubsection{Proof of the global existence}
We now turn to the proof of the global existence. Applying $\triangle_j$ to the first equation of the system (\ref{s2}) yields that
\begin{align*}
(\partial_t+u\cdot \nabla )\triangle_ju+\triangle_j\Pi(u,u)=\triangle_j\mathcal{P}\big((\nabla\cdot\xi)\xi\big)+R_{j1},
\end{align*}
with $R_{j1}=u\cdot \nabla\triangle_ju-\triangle_j(u\cdot \nabla)u\triangleq[u\cdot\nabla,\triangle_j]u.$\\
Taking the $L^2$ inner product of the above equation with $\triangle_{j} u$, we easily get
\begin{align*}
  &\frac{1}{2}\frac{d}{dt}\|\triangle_ju(t)\|_{L^{2}}^2-\frac{1}{2}\int_{\mathbb{R}^2}(\nabla\cdot u)|\triangle_{j}u|^2dx \\ \leq&\|\triangle_{j}u\|_{L^2}\Big(\|\triangle_j\Pi(u,u)\|_{L^{2}}
  +\|\triangle_j\mathcal{P}\big((\nabla\cdot\xi)\xi\big)\|_{L^{2}}+\|R_{j1}\|_{L^{2}}\Big),~j\geq-1.
\end{align*}
Note that $\nabla \cdot u=0$, we get
\begin{align*}
  \|\triangle_ju(t)\|_{L^{2}}\leq\|\triangle_ju_0\|_{L^{2}}+\int_0^t\|\triangle_j\Pi(u,u)\|_{L^{2}}
  +\|\triangle_j\mathcal{P}\big((\nabla\cdot\xi)\xi\big)\|_{L^{2}}+\|R_{j1}\|_{L^{2}}dt',~j\geq-1.
\end{align*}
Multiplying both sides of the above inequality by $2^{js_1}$, taking the $l^{2}$ norm , we obtain
\begin{align}\label{uu'}
  \|u\|_{\widetilde{L}^\infty_t(H^{s_1})}\lesssim
  \|u_0\|_{H^{s_1}}+\|\Pi(u,u)\|_{\widetilde{L}^1_t(H^{s_1})}+
  \|\mathcal{P}\big((\nabla\cdot\xi)\xi\big)\|_{\widetilde{L}^1_t(H^{s_1})}
  +\Big\|2^{js_1}\|R_{j1}\|_{L^1_t(L^{2})}\Big\|_{l^{2}}.
\end{align}
Due to Lemma \ref{jiaohuan}, we get
\begin{align}\label{RR}
  \Big\|2^{js_1}\|R_{j1}\|_{L^1_t(L^{2})}\Big\|_{l^{2}}\lesssim \int_0^t\Big\|2^{js_1}\|R_{j1}\|_{L^{2}}\Big\|_{l^{2}}dt'\lesssim \int_0^t \|\nabla u\|_{L^\infty}\|u\|_{H^{s_1}}dt'.
\end{align}
By virtue of Lemma \ref{pi}, we have
\begin{align}\label{pipi}
  \|\Pi(u,u)\|_{\widetilde{L}^1_t(H^{s_1})}\lesssim  \|\Pi(u,u)\|_{L^1_t(H^{s_1})}\lesssim \int_0^t \|u\|_{C^{0,1}}\|u\|_{H^{s_1}}dt'.
\end{align}
We now focus on the term $\mathcal{P}\big((\nabla\cdot\xi)\xi\big).$ By taking advantage of Bony's decomposition and of Lemmas \ref{T}-\ref{R}, we have
\begin{align}\label{phiphi}
  \|\mathcal{P}\big((\nabla\cdot\xi)\xi\big)\|_{\widetilde{L}^1_t(H^{s_1})}\lesssim&
  \|(\nabla\cdot\xi)\xi\|_{\widetilde{L}^1_t(H^{s_1})}\\\nonumber \lesssim&
  \|\nabla\cdot\xi\|_{L^\infty_t(L^\infty)}\|\xi\|_{\widetilde{L}^1_t(H^{s_1})}
  +\|\xi\|_{L^\infty_t(L^\infty)}\|\nabla\cdot\xi\|_{\widetilde{L}^1_t(H^{s_1})}\\ \nonumber \lesssim&\|\nabla\cdot\xi\|_{L^\infty_t(L^\infty)}\|\xi\|_{\widetilde{L}^1_t(H^{s_2+\frac{3}{2}})}
  +\|\xi\|_{L^\infty_t(L^\infty)}
  \|\nabla\cdot\xi\|_{\widetilde{L}^1_t(H^{s_2+\frac{3}{2}})}\\\nonumber \lesssim&\|\nabla\cdot\xi\|_{L^\infty_t(L^\infty)}
  \|\xi\|_{\widetilde{L}^1_t(H^{s_2+1})}^\frac{3}{4}
  \|\xi\|_{\widetilde{L}^1_t(H^{s_2+3})}^\frac{1}{4}
  +\|\xi\|_{L^\infty_t(L^\infty)}
  \|\xi\|_{\widetilde{L}^1_t(H^{s_2+1})}^\frac{1}{4}
  \|\xi\|_{\widetilde{L}^1_t(H^{s_2+3})}^\frac{3}{4}\\\nonumber
  \lesssim&C_{\sigma}\big(\|\nabla\cdot\xi\|_{L^\infty_t(L^\infty)}^\frac{4}{3}
  +\|\xi\|_{L^\infty_t(L^\infty_t)}^4\big)
  \int_0^t\|\xi\|_{H^{s_2+1}}dt'+
  \sigma\|\xi\|_{\widetilde{L}^1_t(H^{s_2+3})}.
  \end{align}
 Plugging the inequalities (\ref{RR})-(\ref{phiphi}) into (\ref{uu'}), we eventually get
 \begin{align}\label{345}
 \|u\|_{\widetilde{L}^\infty_t(H^{s_1})}\lesssim &\|u_0\|_{H^{s_1}}+\int_0^t\Big(\|u\|_{C^{0,1}}
 \|u\|_{H^{s_1}}\\\nonumber&+C_{\sigma}\big(\|\nabla\cdot\xi\|_{L^\infty_t(L^\infty)}^\frac{4}{3}
  +\|\xi\|_{L^\infty_t(L^\infty_t)}^4\big)
  \|\xi\|_{H^{s_2+1}}\Big)dt'+\sigma
  \|\xi\|_{\widetilde{L}^1_t(H^{s_2+3})}.
  \end{align}

 Similarly, applying $\triangle_j$ to the second equation of the system (\ref{s2})  yields that
\begin{align*}
(\partial_t+u\cdot \nabla -\triangle) \triangle_jz=-\triangle_j\nabla\cdot\big((\nabla\cdot\xi)\xi\big)+R_{j2},
\end{align*}
with $R_{j2}=[u\cdot\nabla,\triangle_j]z,$ where we have used $\nabla \cdot (uz)=u \cdot \nabla z+(\nabla \cdot u)z=u\cdot \nabla  z.$\\
Taking the $L^2$ inner product of the above equation with $\triangle_{j} z$,  we get
\begin{align*}
  &\frac{1}{2}\frac{d}{dt}\|\triangle_jz(t)\|_{L^{2}}^2-\frac{1}{2}\int_{\mathbb{R}^2}(\nabla\cdot u)|\triangle_{j}z|^2dx+\|\nabla\triangle_{j}z\|_{L^2}^2 \\ \leq&\|\triangle_{j}z\|_{L^2}\Big(\|\triangle_j\nabla\cdot\big((\nabla\cdot\xi)\xi\big)\|_{L^{2}}
  +\|R_{j2}\|_{L^{2}}\Big),~j\geq-1.
\end{align*}
Note that $\nabla\cdot u=0,$ $\|\nabla\triangle_{-1}z\|_{L^2}\geq 0,$ and by virtue of Lemma \ref{Bi}, $\|\nabla\triangle_{j}z\|_{L^2}\gtrsim 2^{j}\|\triangle_{j}z\|_{L^2},$ for $j\geq 0.$ Therefore, we have
\begin{align}\label{xingxing}
  \|\triangle_jz(t)\|_{L^{2}}+\int_0^t2^{2j}\|\triangle_jz\|_{L^{2}}dt'\lesssim (1+t)\big(&\|\triangle_jz_0\|_{L^{2}}+\int_0^t
  \|\triangle_j\nabla\cdot\big((\nabla\cdot\xi)\xi\big)\|_{L^{2}}
  \\\nonumber+&\|R_{j2}\|_{L^2}dt'\big),~j\geq-1.
\end{align}
Hence multiplying both sides of the above inequality by $2^{js_2}$ and taking the $l^{2}$ norm, we obtain
\begin{align*}
  &\|z\|_{\widetilde{L}^\infty_t(H^{s_2})}+ \|z\|_{\widetilde{L}^1_t(H^{s_2+2})}\\\lesssim& (1+t)\big(\|z_0\|_{H^{s_2}}+
  \|\nabla\cdot\big((\nabla\cdot\xi)\xi\big)\|_{\widetilde{L}^1_t(H^{s_2})}
  +\Big\|2^{js_2}\|R_{j2}\|_{L^1_t(L^{2})}\Big\|_{l^{2}}\big).
\end{align*}
In view of Lemma \ref{jiaohuan}, we get
\begin{align}\label{R2}
   \Big\|2^{js_2}\|R_{j2}\|_{L^1_t(L^{2})}\Big\|_{l^{2}}\lesssim
    &\int_0^t\Big\|2^{js_2}\|R_{j2}\|_{L^{2}}\Big\|_{l^{2}}dt')\\\nonumber\lesssim &\int_0^t\Big(\|\nabla u\|_{L^\infty}\|z\|_{H^{s_2}}+\|\nabla z\|_{L^\infty}\|\nabla u\|_{H^{s_2-1}}\Big)dt'\\\nonumber
   \lesssim &\int_0^t\Big(\|\nabla u\|_{L^\infty}\|z\|_{H^{s_2}}+\|\nabla z\|_{L^\infty}\| u\|_{H^{s_1}}\Big)dt'.
\end{align}
According to Lemmas \ref{s_2+1}, we have
\begin{align}\label{n2}
  &\|\nabla\cdot\big((\nabla\cdot\xi)\xi\big)\|_{\widetilde{L}^1_t(H^{s_2})}
  \lesssim\|(\nabla\cdot\xi)\xi\|_{\widetilde{L}^1_t(H^{s_2+1})}\\\nonumber
  \lesssim&\|\nabla\cdot\xi\|_{L^\infty_t(L^\infty)}\|\xi\|_{\widetilde{L}^1_t(H^{s_2+1})}
  +\|\xi\|_{L^\infty_t(L^\infty)}\|\nabla\cdot\xi\|_{\widetilde{L}^1_t(H^{s_2+1})},\\\nonumber
  \lesssim&\|\nabla\cdot\xi\|_{L^\infty_t(L^\infty)}\int_0^t\|\xi\|_{H^{s_2+1}}dt'
  +C_{\sigma}\|\xi\|_{L^\infty_t(L^\infty)}^2\int_0^t\|\xi\|_{H^{s_2+1}}dt'+\sigma\|\xi\|_{\widetilde{L}^1_t(H^{s_2+3})}
.
\end{align}
Inserting the inequalities (\ref{R2})-(\ref{n2}) into (\ref{xingxing}), we finally get
\begin{align}\label{567}
 &\|z\|_{\widetilde{L}^\infty_t(H^{s_2})}+\|z\|_{\widetilde{L}^1_t(H^{s_2+2})}\\\nonumber
\lesssim&(1+t)\Big(\|z_0\|_{H^{s_2}}+\int_0^t\Big(\|\nabla u\|_{L^\infty}\|z\|_{H^{s_2}}+\|\nabla z\|_{L^\infty}\| u\|_{H^{s_1}}\\\nonumber&+\big(\|\nabla\cdot\xi\|_{L^\infty_t(L^\infty)}
  +C_{\sigma}\|\xi\|_{L^\infty_t(L^\infty)}^2\big)\|\xi\|_{H^{s_2+1}}\Big)dt'+\sigma\|\xi\|_{\widetilde{L}^1_t(H^{s_2+3})}
  \end{align}

  To deal with the third equation of the system (\ref{s2}), we have
\begin{align}\label{xi2}
  &\|\mathcal{L}\big(u(\nabla\cdot\xi)\big)\|_{\widetilde{L}^1_t(H^{s_2+1})}\lesssim
  \|u(\nabla\cdot\xi)\|_{\widetilde{L}^1_t(H^{s_2+1})}\\\nonumber
  \lesssim&\|\nabla\cdot\xi\|_{L^\infty_t(L^\infty)}
  \|u\|_{\widetilde{L}^1_t(H^{s_2+1})}+
  \|u\|_{\widetilde{L}^\infty_t(L^q)}\|\nabla\cdot\xi\|_{\widetilde{L}^1_t(B^{s_2+1}_{\frac{2q}{q-2},2})}
  \\\nonumber
  \lesssim&\|\nabla\cdot\xi\|_{L^\infty_t(L^\infty)}
  \|u\|_{\widetilde{L}^1_t(H^{s_1})}+
  \|u\|_{L^\infty_t(L^q)}
  \|\nabla\cdot\xi\|_{\widetilde{L}^1_t(H^{s_2+1+\frac{2}{q}})}
  \\\nonumber
  \lesssim&\|\nabla\cdot\xi\|_{L^\infty_t(L^\infty)}
  \int_0^t\|u\|_{H^{s_1}}dt'+
  C_{\sigma}\|u\|_{L^\infty_t(L^q)}^{\frac{2q}{q-2}}
  \int_0^t\|\xi\|_{H^{s_2+1}}dt'+\sigma\|\xi\|_{\widetilde{L}^1_t(H^{s_2+3})},
\end{align}
with $2\leq q < \infty.$
\begin{align}\label{xi3}
  &\|\mathcal{L}(z\xi)\|_{\widetilde{L}^1_t(H^{s_2+1})}\lesssim
  \|z\xi\|_{\widetilde{L}^1_t(H^{s_2+1})}\\\nonumber
  \lesssim&\|\xi\|_{L^\infty_t(L^\infty)}
  \|z\|_{\widetilde{L}^1_t(H^{s_2+1})}+
  \|z\|_{L^\infty_t(L^\infty)}\|\xi\|_{\widetilde{L}^1_t(H^{s_2+1})}
  \\\nonumber
  \lesssim&C_{\sigma}\|\xi\|_{L^\infty_t(L^\infty)}^2
  \int_0^t\|z\|_{H^{s_2}}dt'+\sigma\|z\|_{\widetilde{L}^1_t(H^{s_2+2})}
+
  \|z\|_{L^\infty_t(L^\infty)}\int_0^t\|\xi\|_{H^{s_2+1}}dt'.
\end{align}
Hence,
   \begin{align}\label{568}
  &\|\xi\|_{\widetilde{L}^\infty_t(H^{s_2+1})}+\|\xi\|_{\widetilde{L}^1_t(H^{s_2+3})}\\\nonumber
  \lesssim&(1+t)\Big(\|\xi_0\|_{H^{s_2+1}}
  +\|\mathcal{L}\big(u(\nabla\cdot\xi)\big)\|_{\widetilde{L}^1_t(H^{s_2+1})}
  +\|\mathcal{L}(z\xi)\|_{\widetilde{L}^1_t(H^{s_2+1})}\Big)
  \\\nonumber
\lesssim&(1+t)\Big(\|\xi_0\|_{H^{s_2+1}}+\big(\|\nabla\cdot\xi\|_{L^\infty_t(L^\infty)}
  +
  C_{\sigma}\|u\|_{L^\infty_t(L^q)}^{\frac{2q}{q-2}}+C_{\sigma}\|\xi\|_{L^\infty_t(L^\infty)}^2
  +\|z\|_{L^\infty_t(L^\infty)}\big)\\\nonumber&\times\int_0^t\big(\|u\|_{H^{s_1}}+\|\xi\|_{H^{s_2+1}}+\|z\|_{H^{s_2}}\big)dt'
  +\sigma\big(\|\xi\|_{\widetilde{L}^1_t(H^{s_2+3})}+\|z\|_{\widetilde{L}^1_t(H^{s_2+2})}\big) \Big).
  \end{align}
  Combining (\ref{567}), (\ref{568}) and (\ref{345}), we get $\forall t\in[0,T^*),$
  \begin{align}
    &\|u\|_{\widetilde{L}^\infty_t(H^{s_1})}+\|z\|_{\widetilde{L}^\infty_t(H^{s_2})\cap \widetilde{L}^1_t(H^{s_2+2})}+\|\xi\|_{\widetilde{L}^\infty_t(H^{s_2+1})\cap \widetilde{L}^1_t(H^{s_2+3})}\\
    \nonumber\lesssim
    &(1+t)\big(\|u_0\|_{H^{s_1}}+\|z_0\|_{H^{s_2}}+\|\xi_0\|_{H^{s_2+1}}\big)+(1+t)\int_0^t\Big(
    \|u\|_{C^{0,1}}+\| \nabla z\|_{L^\infty}+\\
    \nonumber&C_{\sigma}\big(\|\nabla\cdot\xi\|_{L^\infty_t(L^\infty)}^\frac{4}{3}
  +\|\xi\|_{L^\infty_t(L^\infty_t)}^4+\|\nabla\cdot\xi\|_{L^\infty_t(L^\infty)}
  +\|\xi\|_{L^\infty_t(L^\infty)}^2+\|\nabla\cdot\xi\|_{L^\infty_t(L^\infty)}
  \\\nonumber&+
  \|u\|_{L^\infty(L^q)}^{\frac{2q}{q-2}}+\|\xi\|_{L^\infty_t(L^\infty)}^2
  +\|z\|_{L^\infty_t(L^\infty)}\big)\Big)
   \times\big(\|u\|_{H^{s_1}}+\|z\|_{H^{s_2}}+\|\xi\|_{H^{s_2+1}}\big)dt'
   \\\nonumber&+(1+t)\sigma\big(\|z\|_{\widetilde{L}^1_t(H^{s_2+2})}+\|\xi\|_{\widetilde{L}^1_t(H^{s_2+3})}\big).
  \end{align}
  Choose $\sigma=c(1+T^*)^{-1}.$
  Lemmas \ref{1l}-\ref{6l} and the inequality (\ref{uq}) imply that
  \begin{align*}
    &C_{\sigma}\big(\|\nabla\cdot\xi\|_{L^\infty_t(L^\infty)}^\frac{4}{3}
  +\|\xi\|_{L^\infty_t(L^\infty_t)}^4+\|\nabla\cdot\xi\|_{L^\infty_t(L^\infty)}
  +\|\xi\|_{L^\infty_t(L^\infty)}^2+\|\nabla\cdot\xi\|_{L^\infty_t(L^\infty)}
  \\\nonumber&+
  \|u\|_{L^\infty(L^q)}^{\frac{2q}{q-2}}+\|\xi\|_{L^\infty_t(L^\infty)}^2
  +\|z\|_{L^\infty_t(L^\infty)}\big)\leq C(T^*)<\infty,
  \end{align*}
  from which it follows that
  \begin{align}
    &\|u\|_{\widetilde{L}^\infty_t(H^{s_1})}+\|z\|_{\widetilde{L}^\infty_t(H^{s_2})}
    +\|\xi\|_{\widetilde{L}^\infty_t(H^{s_2+1})}\\
    \nonumber\leq
    &C(T^*)\Big(\big(\|u_0\|_{H^{s_1}}+\|z_0\|_{H^{s_2}}+\|\xi_0\|_{H^{s_2+1}}\big)+\int_0^t\big(
    \|u\|_{C^{0,1}}+\| \nabla z\|_{L^\infty}+1)\\\nonumber&~~~~~~
   \times\big(\|u\|_{\widetilde{L}^\infty_{t'}(H^{s_1})}+\|z\|_{\widetilde{L}^\infty_{t'}(H^{s_2})}
    +\|\xi\|_{\widetilde{L}^\infty_{t'}(H^{s_2+1})}\big)dt'\Big)\\\nonumber\triangleq &B(t).
  \end{align}
Denote
  \begin{align*}
  &\|u_0\|_{H^{s_1}}+\|z_0\|_{H^{s_2}}+
    \|\xi_0\|_{H^{s_2+1} }\triangleq A_0,\\
    &\|u\|_{\widetilde{L}^\infty_T(H^{s_1})}+\|z\|_{\widetilde{L}^\infty_T(H^{s_2})}+
    \|\xi\|_{\widetilde{L}^\infty_T(H^{s_2+1}) }\triangleq A(t),\end{align*}
Let $\epsilon =\min(1,s_1-2)$ and $\Gamma(r)=1+\log r:[1,\infty)\rightarrow[0,\infty)$
be the function associated with the modulus of continuity $\mu(r) = r(1-\log r).$ We can extend the domain of
definition of $\Gamma$ to $[0,\infty)$ with $\Gamma(s)=\Gamma(1)=1,$ for $0\leq s<1.$
 The function $G(y) \overset{def} {=}\int_1^y\frac{dy'}{\Gamma(y'^{\frac{1}{\varepsilon}})y'}=\varepsilon\log(1+\frac{1}{\varepsilon}\log y)$
then maps $[1,+\infty)$ onto and one-to-one
$[0,+\infty).$\\
Assuming that $A_0>0,$ otherwise $(0,0,0)$ is the global solution. Using Lemma \ref{LL} with $\Lambda = A_0$, we get
\begin{align*}
  B(t)&\leq  C(T^*)\Big(A_0+\int_0^t(\|u\|_{L^\infty}+\|\nabla u\|_{L^\infty}+\| \nabla z\|_{L^\infty}+1)B(t')dt'\Big) \\ &\leq C(T^*)\Big(A_0+\int_0^t\Big\{\|u\|_{B^1_{\infty,\infty}}+\| \nabla z\|_{L^\infty}+1+C_{\epsilon}\Big(\|u\|_{C_{\mu}}+A_0\Big)\Big(1+
  \Gamma\Big(\big(\frac{\|\nabla u\|_{C^{0,\epsilon}}}{\|u\|_{C_{\mu}}+A_0}\big)
  ^{\frac{1}{\epsilon}}\Big)\Big)\Big\}B(t')dt'\Big)\\ &\leq C(\epsilon,T^*)\Big(A_0+\int_0^t\Big\{\big(\|u\|_{B^1_{\infty,\infty}}+\| \nabla z\|_{L^\infty}+1+A_0\big)\Big(1+
  \Gamma\Big(\big(\frac{C\| u(t')\|_{H^{s_1}}}{A_0}\big)
  ^{\frac{1}{\epsilon}}\Big)\Big)\Big\}B(t')dt'\Big)\\ &\leq C(\epsilon,T^*)\Big(A_0+\int_0^t\Big(\|u\|_{B^1_{\infty,\infty}}+\| \nabla z\|_{L^\infty}+1+A_0\Big)
  \Gamma\Big(\big(\frac{CB(t')}{A_0}\big)
  ^{\frac{1}{\epsilon}}\Big)\Big)B(t')dt'\Big)\\
  &\triangleq\frac{R(t)A_0}{C},
\end{align*}
 where we have used $B^1_{\infty,\infty}\hookrightarrow L^\infty,$ $B^1_{\infty,\infty}\hookrightarrow C_{\mu},$ $H^{s_1}\hookrightarrow C^{0,\epsilon},$ $\| u(t')\|_{H^{s_1}}\leq B(t')$ and $C$ has been chosen large enough such that $R(t)=\frac{B(t)C}{A_0}\geq C>1.$\\
 Because the function $\Gamma$ is nondecreasing, after a few computations, we have that
\begin{align*}
  \frac{d}{dt}R(t)
  \leq&\Gamma\big(R(t)
  ^{\frac{1}{\epsilon}}\big)R(t) C(\epsilon,T^*)\big(\|u(t)\|_{B^1_{\infty,\infty}}+\| \nabla z(t)\|_{L^\infty}+1+A_0\big),
  \end{align*}
thus $$\frac{d}{dt}G(R(t))\leq C(\epsilon,T^*)(\|u(t)\|_{B^1_{\infty,\infty}}+A_0+1+\|\nabla z(t)\|_{L^\infty}).$$
 Integrating then gives
 \begin{align*}
   R(t)\leq G^{-1}\Big(G\big(R(0)\big)+\int_0^tC(\epsilon,T^*)(\|u\|_{B^1_{\infty,\infty}}+A_0+1+\|\nabla z\|_{L^\infty})dt'\Big)<\infty,
 \end{align*}
 where we have used Lemmas \ref{5l}-\ref{6l}.
Therefore, $
  \|u(t)\|_{H^{s_1}},$ $\|z(t)\|_{H^{s_2}},$ and $
    \|\xi(t)\|_{H^{s_2+1} }$ stay bounded on $[0,T^*)$. The local existence part of
  Theorem \ref{a1} then enables us to extend the solution beyond $T^*,$ which stands in contradiction to the definition of $T^*.$ Hence $T^*=+\infty.$
  This completes the proof of the theorem.\qed

 \section{Proof of Theorem \ref{a3}}
To begin, we denote by $BMO$ the space of functions of bounded mean oscillations. It is well known that $BMO$ strictly includes $L^\infty$. We introduce the following Hardy-Littlewood-Sobolev inequality.
\begin{lemm}\cite{Lemari¨¦-Rieusset}\label{Lemari¨¦-Rieusset}
For $0<\gamma<d,$ the operator $(-\triangle)^{\frac{\gamma}{2}}$ is bounded from the Hardy space $\mathcal{H}^1$ to $L^{\frac{d}{d-\gamma}}$ and from $L^{\frac{d}{\gamma}}$ to $BMO.$
\end{lemm}
 \subsection{ Global existence for the $ENPP$ system}
Let $(u_0,z_0,\xi_0)=\big(u_0,n_0+p_0,-\nabla(\triangle)^{-1}(n_0-p_0)\big).$ According to Theorem \ref{a1}, there exists a global solution $(u,z,\xi)$ satisfies the system (\ref{s2}) in the spaces defined as in Theorem \ref{a1}. Since $\mathcal{L}\xi=\xi,$ it is then easy to that $(u,n,p)=(u,\frac{z+\nabla\cdot \xi}{2},\frac{z-\nabla\cdot \xi}{2})$ solves the system (\ref{s3}).\\
Denote $$\phi_0\triangleq-(-\triangle)^{-1}\nabla \cdot \xi.$$ As $\xi\in L^\infty(\mathbb{R}^+;L^2),$ applying Lemma \ref{Lemari¨¦-Rieusset} with $d=2$ and $\gamma=1$ implies that $\phi_0\in L^\infty(\mathbb{R}^+;BMO).$ Thanks again to the fact that $\mathcal{L}\xi=\xi,$ we have $\nabla \phi_0=\mathcal{L}\xi=\xi,$ and $\triangle \phi_0=\nabla \cdot \xi=n-p.$
Similarly, let $$P_0\triangleq P_{\pi}(u,u)-(-\triangle)^{-1}\nabla \cdot\big((\nabla \cdot \xi)\xi\big),$$
  where $P_{\pi}(u,u)\in L^\infty(\mathbb{R}^+;H^{s_1+1})$ is defined as in Lemma \ref{pi}.
 Note that $\xi\in L^\infty(\mathbb{R}^+;H^{s_1+1})$ with $s_1>1$ implies that $\nabla \cdot \xi\in L^\infty(\mathbb{R}^+;L^2)$ and $\xi\in L^\infty(\mathbb{R}^+;L^\infty)$. Again using lemma \ref{Lemari¨¦-Rieusset}, we get
 $$P_0\in L^\infty(\mathbb{R}^+;H^{s_1+1}+BMO)\hookrightarrow L^\infty(\mathbb{R}^+;L^\infty+BMO)\hookrightarrow L^\infty(\mathbb{R}^+;BMO).$$
 Finally, it is easy to see that $(u,\frac{z+\nabla\cdot \xi}{2},\frac{z-\nabla\cdot \xi}{2},P_0,\phi_0)$ satisfies the $ENPP$ system.
 \subsection{ Uniqueness for the $ENPP$ system}
 Suppose that there exists a global solution $(u,n,p,P,\phi)$ satisfing the $ENPP$ system in the spaces defined as in Theorem \ref{a3}. We first show that $$\nabla \Phi=-\nabla(-\triangle)^{-1}(n-p)\triangleq\xi,~\textit{and} ~\nabla P=\pi(u,u)+(I-\mathcal{P})\big((n-p)\nabla(-\triangle)^{-1}(p-n)\big).$$
 In fact, Let $\phi_0,~P_0$ be defined as in the above subsection. As $\triangle \phi=n-p=\triangle \phi_0,$ hence $\phi-\phi_{0}$ is a harmonic polynomial. Note that $\phi \in L^\infty(\mathbb{R}^+;BMO)$ is required in Theorem \ref{a3} and $\phi_0 \in L^\infty(\mathbb{R}^+;BMO)$ is illustrated before. Thus $\phi-\phi_0$ depends only on t, and
  \begin{align}\label{pphi}
    \nabla \phi=\nabla\phi_{0}=\xi=-\nabla(-\triangle)^{-1}(n-p).
  \end{align}
  Next applying the operator $\nabla\cdot$ to the first equation of the $ENPP$ system, we get
  \begin{align*}
    -\triangle P=\nabla\cdot(u\cdot\nabla u)-\nabla\cdot((\nabla \cdot \xi)\xi)
    =-\triangle P_0.
  \end{align*}
  Note that $P-P_0$ is in
  $L^\infty(\mathbb{R}^+;BMO).$ Similar arguments as that for $\phi-\phi_0$ yield that
  \begin{align*}
    \nabla P=\nabla P_0
    =\Pi(u,u)-\nabla(-\triangle)^{-1}\nabla \cdot\big((\nabla \cdot \xi)\xi\big)=\Pi(u,u)+(I-\mathcal{P})\big((n-p)\nabla(-\triangle)^{-1}(p-n)\big).
  \end{align*}
Next it is easy to see that $(u,n,p,\xi)$ solves the system (\ref{s3}), and $(u,n+p,\xi)$ solves the system (\ref{s2}). The uniqueness of the system (\ref{s2}) in Theorem \ref{a1} then implies that $(u,n,p,\nabla P,\nabla \phi)$ is uniquely determined by the initial data.
This completes the proof of the theorem. \qed\\
%\bigskip
\noindent\textbf{Acknowledgements}. This work was partially supported by NNSFC (No. 11271382), RFDP (No. 20120171110014), and
the key project of Sun Yat-sen University.
\phantomsection
\addcontentsline{toc}{section}{\refname}
%\bibliographystyle{unsrtnat}
%\nocite{*}
\bibliography{reference}

\end{document}